%% file: bhps.tex
\documentclass[reqno,oneside,12pt]{amsart}

\usepackage[T1]{fontenc}
\usepackage{times,mathptm}
\usepackage{amssymb,epsfig,verbatim,xypic}
\theoremstyle{plain}

\newtheorem{thm}{Theorem}[section]

\newtheorem{cor}[thm]{Corollary}
\newtheorem{pro}[thm]{Proposition}
\newtheorem{lem}[thm]{Lemma}
\newtheorem{proposition-principale}[thm]{Proposition principale}
\newtheorem{thm-principal}{Th\'eor\`eme principal}[section]

\theoremstyle{definition}

\newtheorem{eg}[thm]{Example}
\newtheorem{rem}[thm]{Remark}

\def\vv{\vspace{0.2cm}}

\def\C{\mathbf{C}}
\def\R{\mathbf{R}}

\def\Z{\mathbf{Z}}

\def\P{\mathbb{P}}
\def\Sphere{\mathbb{S}}
\def\disk{\mathbb{D}}

\def\T{\mathbb{T}}

\def\SS{\sf{Fam}}
\def\Sing{{\sf{Sing}}}

\def\Te{{\sf{Teich}}}
\def\Rep{{\sf{Rep}}}

\def\QF{{\sf{QF}}}

\def\Be{{\sf{Bers}}}

\def\disc{\mathbb{D}}
\def\H{\mathbb{H}}

\def\A{\mathcal{A}}
\def\Out{{\sf{Out}}}
\def\Aut{{\sf{Aut}}}

\def\DF{{\sf{DF}}}
\def\Int{{\sf{Int}}}

\def\MCG{{\sf{MCG}}}
\def\PGL{{\sf{PGL}}\,}
\def\PSL{{\sf{PSL}}\,}
\def\GL{{\sf{GL}}\,}
\def\SL{{\sf{SL}}\,}

\def\SU{{\sf{SU}}\,}

\def\tr{{\sf{tr}}}

\def\diam{{\text{diam}}}

\def\Ind{{\text{Ind}}}

\def\gro{{\sf{gro}}}
\def\Gro{{\sf{Gro}}}

\addtocounter{section}{0}             
\numberwithin{equation}{section}       

\begin{document}

\setlength{\baselineskip}{0.51cm}        
%
%
\title[Bers and H\'enon, Painlev\'e and Schr\"odinger]
{Bers and H\'enon, Painlev\'e and Schr\"odinger}
\date{2007}
\author{Serge Cantat}
\address{D\'epartement de math\'ematiques\\
         Universit\'e de Rennes\\
         Rennes\\
         France}
\email{serge.cantat@univ-rennes1.fr}
%
%

%
%

%
%


\begin{abstract} 
In this paper, we pursue the study of the holomorphic dynamics of mapping
class groups on $2$-dimensional character varieties, as initiated in 
\cite{Cantat-Loray:2007, Iwasaki-Uehara:2006}. We shall 
show that the dynamics of pseudo-Anosov mapping classes resembles in many 
ways the dynamics of H\'enon mappings, and then apply this idea to answer
open questions concerning the geometry of discrete and faithful representations, 
Painlev\'e sixth equation, and discrete Schr\"odinger operators.
 \end{abstract}

\maketitle

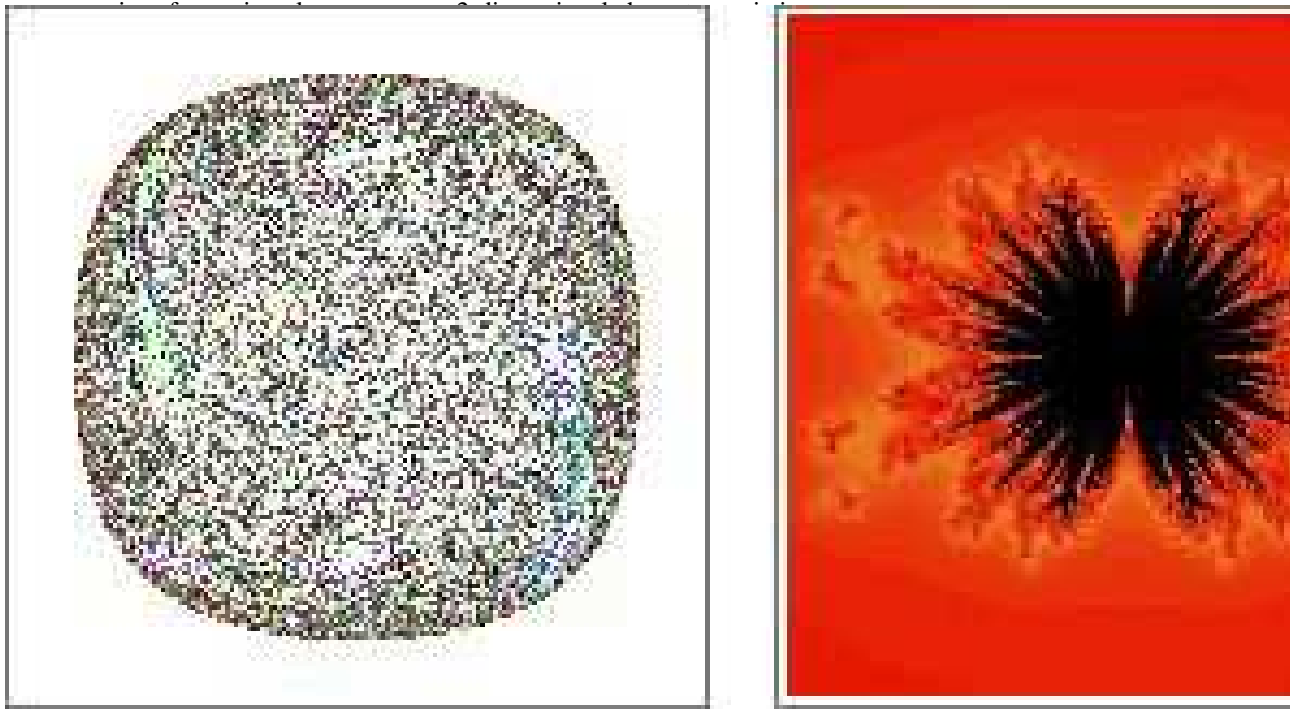
\begin{figure}[h]\label{fig:intro}
\input{intro-julia.pstex_t}
\caption{{\sf{Dynamics on character surfaces.}} {{Left: Dynamics on the real part of a cubic 
surface. Right: A slice of the set of complex points with bounded orbit.}}}
\end{figure}

\setcounter{tocdepth}{1}
\tableofcontents

\section{Introduction}

\subsection{Character variety and dynamics} 

Let $\T_1$ be the once punctured torus. Its fundamental group is isomorphic 
to the free group $F_2=\langle \alpha, \beta \, \vert \, \emptyset \rangle,$ the
commutator of $\alpha$ and $\beta$ corresponding to a simple loop around the 
puncture. Since any representation $\rho:F_2\to \SL(2,\C)$ is uniquely determined by
$\rho(\alpha)$ and $\rho(\beta),$  the set $\Rep(\T_1)$
of representations of $\pi_1(\T_1)$ into $\SL(2,\C)$ is isomorphic to 
$\SL(2,\C)\times \SL(2,\C).$ The group $\SL(2,\C)$ acts on this set by conjugation,
preserving the three traces 
$$
x=\tr(\rho(\alpha)), \quad y=\tr(\rho(\beta)),\quad z=\tr(\rho(\alpha\beta)).
$$ 
It turns out that the map $\chi:\Rep(\T_1)\to \C^3,$ defined
by $\chi(\rho)=(x,y,z),$ realizes an isomorphism between the algebraic
quotient 
$
\Rep(\T_1)/\!/\SL(2,C),
$
where $\SL(2\C)$ acts by conjugation,
and the complex affine space $\C^3.$ 
This quotient will be referred to as the {\sl{character variety
of the once punctured torus}}.

The automorphism group $\Aut(F_2)$ acts by composition on $\Rep(\T_1),$ and induces
an action of the mapping class group 
$$
\MCG^*(\T_1)=\Out(F_2)=\GL(2,\Z)
$$ 
on the character variety $\C^3$ by polynomial diffeomorphisms. Since the 
conjugacy class of the commutator $[\alpha,\beta]$ is invariant under $\Out(F_2),$
this action preserves the level sets of the polynomial function
$\tr(\rho[ \alpha, \beta])=x^2+y^2+z^2 -xyz -2.$
As a consequence, for each complex number $D,$ we get a morphism from 
$\Out(F_2)$ to the group $\Aut(S_D)$ of polynomial diffeomorphisms of the
surface $S_D,$  defined by
$$
x^2+y^2+z^2=xyz+D.
$$
The goal of this paper is to describe the dynamics of mapping classes 
on these surfaces, both on the complex surface $S_D(\C)$ and on the real 
surface $S_D(\R)$ when $D$ is a real number. More generally, {\sl{we shall 
study the dynamics of mapping classes on the character variety of the 
$4$-punctured sphere}}, but we restrict ourselves to the simpler case of the punctured 
torus in the introduction. 


\subsection{H\'enon type dynamics}\label{par:intro-henon}

Let us fix an element $f$ of the mapping class group $\MCG^*(\T_1),$ 
that we view simultaneously as a matrix $M_f$ in $\GL(2,\Z)=\Out(F_2)$ or as a
polynomial automorphism, still denoted $f,$  of the affine space $\chi(\T_1)=\C^3$
 preserving the family of cubic surfaces $S_D.$
Let $\lambda(f)$ be the spectral radius of $M_f,$ so that $f$ is pseudo-Anosov
if and only if $\lambda(f)>1.$  

In \cite{Cantat-Loray:2007} (see also  
\cite{Iwasaki-Uehara:2006}), it is proved that the topological entropy of 
$f:S_D(\C)\to S_D(\C)$ is equal to $\log(\lambda(f))$ for all choices of $D.$ 
The dynamics of
mapping classes with zero entropy is described in details 
in \cite{Goldman:1997,Cantat-Loray:2007}. 
In section \ref{par:Henon}, we shall show that the 
dynamics of pseudo-Anosov classes resembles the dynamics of 
H\'enon automorphisms of the complex plane: All techniques from holomorphic 
dynamics that have been developed for H\'enon automorphisms 
can  be applied to understand the dynamics of mapping classes on 
character varieties !

As a corollary of this principle, we shall get a positive answer to three different
questions. The first one concerns quasi-fuchsian groups and the geometry of
the quasi-fuchsian set. The second one concerns the spectrum of certain discrete 
Schr\"odinger operators, while the third question is related to Painlev\'e sixth equation. 


\subsection{Quasi-Fuchsian spaces and a question of Goldman and Dumas}\label{par:intro-gold}
 First, we answer positively a question of Goldman
and Dumas (see problem 3.5 in \cite{Goldman:survey}),  that we now describe.

When the parameter $D$ is equal to $2,$ the trace of $\rho[\alpha,\beta]$ vanishes,
so that the representations $\rho$ with $\chi(\rho)$ in $S_2(\C)$ send the commutator
$[\alpha,\beta]$ to an element of order $4$ in $\SL(2,\C).$ This means that the surface
$S_2$ corresponds in fact to representations of the group 
$
G=\langle \alpha, \beta\,  \vert\,  [\alpha, \beta]^4\rangle.
$ 
Let $\DF$ be the subset of $S_2(\C)$ corresponding to discrete and faithful representations of $G.$  Some of these representations are
fuchsian: These representations $F_2 \to \SL(2,\R)$ come
from the existence of hyperbolic metrics on $\T_1$ with an orbifold point of 
angle $\pi$ at the puncture. The interior of $\DF$ corresponds to quasi-fuchsian 
deformations of those fuchsian representations.


Let us now consider the set of conjugacy classes of representations 
$\rho:G\to\SU(2).$ This set  coincides with the compact connected component of 
$S_2(\R).$ Typical representations into $\SU(2)$ 
have a dense image and, in this respect,  are quite different from discrete faithful 
representations into $\SL(2,\C).$  The following theorem  shows that orbits of the mapping class 
group may contain both types of representations in their closure. 

\begin{thm}\label{thm:dg} Let $G$ be the finitely presented group 
$\langle \alpha, \beta \, \vert \, [\alpha, \beta]^4\rangle.$
There exists a representation $\rho:G\to \SL(2,\C),$  such that  the closure of
the orbit of its conjugacy class $\chi(\rho)$ under the action of $\Out(F_2)$ 
contains both
\begin{itemize}
\item  the conjugacy class of at least one discrete and faithful 
representation $\rho':G\to\SL(2,\C),$ and 
\item the whole set  of conjugacy classes
of $SU(2)$-representations of the group $G.$
\end{itemize}
\end{thm}

This result answers positively and  precisely the question raised by Dumas and Goldman. 
The strategy of proof is quite general and leads to many other examples; one of 
them is given in \S \ref{par:qfex}. The representation $\rho'$ is very special; it
corresponds to certain discrete representations provided by Thurston's hyperbolization theorem for mapping tori with pseudo-Anosov monodromy. 
The same idea may be used to describe $\DF$ in dynamical terms (see section \ref{par:quasifuchsian}). To sum up, holomorphic dynamics turns out to be useful to
understand the quasi-fuchsian locus and its Bers parametrization.

\subsection{Real dynamics, discrete Schr\"odinger operators, and Painlev\'e VI equation}

The fact that the dynamics of mapping classes is similar to the dynamics of
H\'enon automorphisms will prove useful to study the real dynamics of mapping classes, 
i.e. the dynamics of $f$ on the real part $S_D(\R)$ when $D$ is a real number. 
The following theorem, which  is the main result of section \ref{par:realdynamics}, 
provides a complete answer to a conjecture popularized by Kadanoff some twenty five
years ago (see \cite{Kohmoto-Kadanoff-Tang:1983}, p. 1872, for a somewhat weaker
question).  We refer to papers of Casdagli and Roberts (see \cite{Casdagli:1986} and 
\cite{Roberts:1996}), and references therein for a nice mathematical introduction  
to the subject.

\begin{thm}\label{thm:hypreal}
Let $D$ be a real number.  If $f\in \MCG^*(\T_1)$ is a pseudo-Anosov mapping class,
the topological entropy of $f:S_D(\R)\to S_D(\R)$ is bounded from above 
by $\log(\lambda(f)),$ and the five following properties are equivalent
\begin{itemize}
\item the topological entropy of $f:S_D(\R)\to S_D(\R)$ is equal to $\log(\lambda(f))$;
\item all periodic points of $f:S_D(\C)\to S_D(\C)$ are contained in $S_D(\R)$;
\item the topological entropy of $f:S_D(\R)\to S_D(\R)$ is positive and the
dynamics of $f$ on the set 
$K(f,\R)=\{m \in S_D(\R)\, \vert \, (f^n(m))_{n\in \Z} \, {\text{ is bounded }}\}$
is uniformly hyperbolic;
\item the surface $S_D(\R)$ is connected
\item the real parameter $D$ is greater than or equal to $4.$
\end{itemize}
\end{thm}

The main point is the fact that the dynamics is uniformly hyperbolic when $D\geq 4.$ 
As we shall explain in section \ref{par:schrodinger}, this may be used to study the 
spectrum of discrete Schr\"odinger operators, the potential of which 
is generated by a primitive substitution: For example, we shall show that the 
Hausdorff dimension of the spectrum  is positive but strictly less than $1$ 
(see section \ref{par:schrodinger} for precise results). 
This gives also examples of Painlev\'e VI equations with nice and rich 
monodromy (see section \ref{par:painleve6}), thereby answering a question of
Iwasaki and Uehara in \cite{Iwasaki:NagoyaQuestion}.

\begin{table}[htdp]
\caption{Dynamics of pseudo-Anosov classes on $S_D(\R)$}
\vspace{-0.2cm}
\begin{center}
\begin{tabular}{|c|c|c|}
\hline 
values of parameter & real part of $S_D$ & dynamics on $K(f,\R)$ \\
\hline
\hline
\vspace{-0.1cm} & \vspace{-0.1cm} & \vspace{-0.1cm} \\
$D<0$ & four disks  &  $K(f, \R)=\emptyset$ \\
\vspace{-0.1cm} & \vspace{-0.1cm} & \vspace{-0.1cm} \\
\hline 
\vspace{-0.1cm} & \vspace{-0.1cm} & \vspace{-0.1cm} \\
$D=0$ & four disks and a point &  $K(f, \R)=\{(0,0,0)\}$ \\
\vspace{-0.1cm} & \vspace{-0.1cm} & \vspace{-0.1cm} \\
\hline
\vspace{-0.1cm} & \vspace{-0.1cm} & \vspace{-0.1cm} \\
$0< D < 4$ &  $\;$four disks and a sphere$\;$   & non uniformly hyperbolic\\
\vspace{-0.1cm} & \vspace{-0.1cm} & \vspace{-0.1cm} \\
\hline 
\vspace{-0.1cm} & \vspace{-0.1cm} & \vspace{-0.1cm} \\
$D=4$ &  the Cayley cubic & uniformly hyperbolic \\
\vspace{-0.1cm} & \vspace{-0.1cm} & \vspace{-0.1cm} \\
\hline 
\vspace{-0.1cm} & \vspace{-0.1cm} & \vspace{-0.1cm} \\
$D>4$  & a connected surface  & uniformly hyperbolic \\ 
\vspace{0.0cm} & \vspace{0.0cm} & \vspace{0.0cm} \\
\hline 
\end{tabular}
\end{center}
\end{table}


\subsection{Organization of the paper}  
As mentioned above, we shall study the dynamics of the mapping class group
of the four punctured sphere on its character variety; this includes the case of 
the once punctured torus as a particular case.
Section \ref{par:prelim} summarizes known useful results, fixes the notations, 
and describes the dynamics of mapping classes at infinity. 
Section \ref{par:Henon} establishes a dictionary between the H\'enon case
and the case of character varieties, listing important consequences regarding the dynamics
of mapping classes. This is applied in section \ref{par:quasifuchsian} to study the quasi-fuchsian space. 
Section \ref{par:realdynamics} describes the  dynamics
of mapping classes on the real algebraic surfaces $S_D(\R),$ for $D\in \R.$  This is certainly 
the most involved part of this paper. It requires a translation of most known facts for H\'enon 
automorphisms to the case of  character varieties, and a study of one parameter
families of real polynomial automorphisms with maximal entropy. The proof of
theorem \ref{thm:hypreal}, which is given in sections \ref{par:rdmeqh} and 
\ref{par:rduh}, could also be used in the study of families of H\'enon mappings. 
We then apply theorem \ref{thm:hypreal} to the study of Schr\"odinger operators  
and  Painlev\'e VI equations in section \ref{par:schrodinger}.


\subsection{Acknowledgement} This paper greatly benefited from 
discussions with Frank Loray, with whom I collaborated on a closely
related article (see \cite{Cantat-Loray:2007}). I also want to thank Eric Bedford, 
Cliff Earle, Bill Goldman, Katsunori Iwasaki, Robert MacKay, Yair Minsky, John Smillie, 
Takato Uehara and Karen Vogtmann for illuminating talks and useful discussions. 
Most of the content of this paper has been written while I was visiting Cornell 
University in 2006/2007, and part of it was already described during a conference of the 
ACI "Syst\`emes Dynamiques Polynomiaux" in 2004: I thank both institutions for
their support.


\section{The character variety of the four punctured sphere and its automorphisms}\label{par:prelim}

This section summarizes known results concerning the character variety of a
four punctured sphere and the action of its mapping class group 
on this algebraic variety. Most of these results can be found in 
\cite{Benedetto-Goldman:1999}, \cite{Iwasaki-Uehara:2006}, and \cite{Cantat-Loray:2007}.


\subsection{The sphere minus four points}\label{par:prelim-equations}

Let $\Sphere^2_4$ be the four punctured sphere. Its fundamental group
is isomorphic to a free group of rank $3,$
$$
\pi_1(\Sphere^2_4)= \langle \alpha,\beta,\gamma,\delta \, \vert \, \alpha\beta\gamma\delta = 1 \rangle,
$$
where the four homotopy classes $\alpha,$ $\beta,$ $\gamma,$ and $\delta$ 
correspond to loops around the puncture. 
Let $\Rep(\Sphere^2_4)$ be the set of representations of $\pi_1(\Sphere^2_4)$ into $\SL(2,\C).$ Let us associate the $7$ following traces to any element $\rho$ of $\Rep(\Sphere^2_4),$
$$
\begin{array}{clclclc}
 a = \tr(\rho(\alpha)) & ; & b  =  \tr(\rho(\beta)) & ; & c = \tr(\rho(\gamma)) & ; &  d  =  \tr(\rho(\delta)) \\
 x  = \tr(\rho(\alpha \beta)) & ; &  y  =  \tr(\rho(\beta \gamma))   & ; & z  =  \tr( \rho(\gamma \alpha)) . & \, & \,
 \end{array}
$$
The polynomial map $\chi:\Rep(\Sphere^2_4)\to \C^7$ defined by 
$\chi(\rho)=(a,b,c,d,x,y,z)$
is invariant under conjugation, by which we mean that $\chi(\rho')=\chi(\rho)$ 
if $\rho'$ is conjugate to $\rho$ by an element of $\SL(2,\C),$ and it turns out 
that  the algebra of polynomial functions on $\Rep(\Sphere^2_4)$ which are invariant
under conjugation is generated by the components of $\chi.$
Moreover, the components of $\chi$ satisfy the quartic equation 
\begin{equation}\label{eq:surface}
x^2+y^2+z^2+x y z = A x + By +C z + D,
\end{equation}
in which the variables $A,$ $B,$ $C,$ and $D$ are given by
\begin{equation}\label{eq:parameters}
\begin{array}{c} A  =  ab + cd, \quad B  =  ad + bc,\quad  C  =  ac + bd,  \\
 {\text{and}} \quad  D  =  4 - a^2 - b^2- c^2 - d^2 - abcd .  \end{array} 
\end{equation}
In other words,  the algebraic quotient 
$
\chi(\Sphere^2_4):=\Rep(\Sphere^2_4)/\! /\SL(2,\C)
$ 
of  $\Rep(\Sphere^2_4)$ by the action of $\SL(2,\C)$ by conjugation
 is isomorphic to the six-dimensional quartic hypersurface of $\C^7$ 
 defined by equation (\ref{eq:surface}). 

The affine algebraic variety   $\chi(\Sphere^2_4)$ is called the ``character
variety of $\Sphere^2_4$''. For each choice of four complex parameters $A,$ 
$B,$ $C,$ and $D,$ $S_{(A,B,C,D)}$ (or $S$ is there is no obvious possible 
confusion) will denote the cubic surface of $\C^3$ defined by the equation 
(\ref{eq:surface}). The family of surfaces  $S_{(A,B,C,D)},$ with
$A,$ $B,$ $C,$ and $D$ describing $\C,$ will be denoted by $\SS.$


\subsection{Automorphisms and the modular group $\Gamma_2^*$}\label{par:prelim-amg}

The (extended) mapping class group of $\Sphere^2_4$ acts on $\chi(\Sphere^2_4)$ 
by polynomial automorphisms: This defines a morphism 
$$
\left\{ \begin{array}{ccc} \Out(\pi_1(\Sphere^2_4)) & \to & \Aut(\chi(\Sphere^2_4)) \\
 \Phi & \mapsto & f_\Phi
\end{array}\right.
$$
such that $f_\Phi ( \chi (\rho) ) = \chi (\rho\circ \Phi^{-1})$ for any representation $\rho.$

The group $\Out(\pi_1(\Sphere^2_4))$ contains a copy of $\PGL(2,\Z)$ which is obtained
as follows. Let $\T=\R^2/\Z^2$ be a torus and $\sigma$ be the
involution of $\T$ defined by $\sigma(x,y)=(-x,-y).$ The fixed point set of $\sigma$ is 
the $2$-torsion subgroup of $\T.$ The quotient $\T/\sigma$ is homeomorphic to the 
sphere, $\Sphere^2,$ and the quotient map $\pi:\T\to \T/\sigma=\Sphere^2$ 
 has four ramification points, corresponding to the four fixed points of $\sigma.$ 
The group $\GL(2,\Z)$ acts linearly on $\T$ and commutes with $\sigma.$ This 
yelds an action of $\PGL(2,\Z)$  on the sphere $\Sphere^2,$ which permutes the ramification points 
of $\pi.$ Taking these four ramification points as the punctures of $\Sphere^2_4,$we 
get a morphism 
$$
\PGL(2,\Z)\to \MCG^*(\Sphere^2_4),
$$
that turns out to be injective, with finite index image (see \cite{Birman:Book,Cantat-Loray:2007}). 
As a consequence, $\PGL(2,\Z)$ acts by polynomial 
transformations on $\chi(\Sphere^2_4).$ 

Let $\Gamma_2^*$ be the subgroup of $\PGL(2,\Z)$ whose elements coincide 
with the identity  modulo $2.$ This group coincides with the 
stabilizer of the fixed points of $\sigma,$ so that $\Gamma_2^*$ acts 
on $\Sphere^2_4$ and fixes its four punctures. Consequently, $\Gamma_2^*$ 
acts polynomially on $\chi(\Sphere^2_4)$ and preserves the fibers of the projection
$$
(a,b,c,d,x,y,z)\mapsto (a,b,c,d).
$$
From this we obtain, for any choice of four complex parameters $(A,B,C,D),$ a morphism from 
$\Gamma_2^*$ to the group $\Aut(S_{(A,B,C,D)})$ of polynomial diffeomorphisms of 
the surface $S_{(A,B,C,D)}.$

\begin{thm}[\cite{El-Huti:1974,Cantat-Loray:2007}]
For any choice of  $A,$ $B,$ $C,$ and $D,$ the morphism 
$$
\Gamma_2^*\to \Aut(S_{(A,B,C,D)})  
$$
is injective and the index of its image is bounded from above by $24.$ For a generic 
choice of the parameters, this morphism is an isomorphism. 
\end{thm}

The area form $\Omega,$ which is globally defined by the formulas 
$$
\Omega = \frac{dx\wedge dy}{2z+xy-C} = \frac{dy\wedge dz}{2x + yz -A} 
= \frac{dz\wedge dx}{2y+zx-B}
$$
on $S\setminus\Sing(S),$ is almost invariant under the action of $\Gamma_2^*,$ 
by which we mean that $f^*\Omega = \pm \Omega$  for any $f$ in $\Gamma_2^*$ 
(see \cite{Cantat-Loray:2007}). In particular, the dynamics of 
mapping classes on each surface $S$ is conservative.


\subsection{Compactification and automorphisms}\label{par:prelim-ca}
Let $S$ be any member of the family $\SS.$ The closure $\overline{S}$ 
of $S$ in $\P^3(\C)$ is given by the cubic homogeneous equation 
$$
w(x^2+y^2+z^2)+ xyz = w^2(Ax+By+Cz) + Dw^3.
$$
As a consequence, one easily proves that the trace of $\overline{S}$ 
at infinity does not depend on the parameters and coincides with the 
triangle $\Delta$ given by the equations
$$
xyz=0, \quad w=0,
$$
and, moreover, that  the surface $\overline{S}$ is smooth in a neighborhood 
of $\Delta$ (all singularities of $\overline{S}$ are contained in $S$).
The three sides of $\Delta$ are the lines $D_x=\{x=0,w=0\},$
$D_y=\{y=0,w=0\}$ and $D_z=\{z=0,w=0\}$; the vertices are $v_x=[1:0:0:0],$
$v_y=[0:1:0:0]$ and $v_z=[0:0:1:0].$ 
The ``middle points'' of the sides are respectively
$m_x=[0:1:1:0],$ $m_y=[1:0:1:0],$ and $ m_z=[1:1:0:0].$

Since the equation defining $S$ is of degree $2$ with respect to the $x$ variable, 
each point $(x,y,z)$ of $S$ gives rise to a unique second point $(x',y,z).$ This procedure
determines a holomorphic involution of $S,$ namely
$s_x(x,y,z)=(A-yz-x,y,z).$
Geometrically, the involution $s_x$ corresponds to the following: If $m$ is a point of
$\overline{S},$ the projective line which joins $m$ and the vertex $v_x$ of the 
triangle $\Delta$ intersects ${\overline{S}}$ on a third point; this point is $s_x(m).$ 
The same construction provides two more involutions $s_y$ and $s_z,$ and therefore  
a subgroup 
$$
\A = \langle s_x,s_y,s_z\rangle
$$
of the group $\Aut(S)$ of polynomial automorphisms of the surface $S.$ 
It is proved in \cite{Cantat-Loray:2007} that the group $\A$ coincides with the 
image of $\Gamma_2^*$ into $\Aut(S),$ that is obtained by the action of 
$\Gamma_2^*\subset \MCG^*(\Sphere^2_4)$ on the character variety $\chi(\Sphere^2_4).$ More 
precisely, $s_x,$ $s_y,$ and $s_z$ correspond respectively to the automorphisms 
determined by the following elements of $\Gamma^*_2$
$$
r_x=\left(\begin{array}{cc} -1 & -2 \\ 0 & 1 \end{array}\right), \quad r_y= 
\left(\begin{array}{cc} -1 & 0 \\ 0 & 1 \end{array}\right), \quad r_z= 
\left(\begin{array}{cc} 1 & 0 \\ -2 & -1 \end{array}\right).
$$  
In particular, there is no non trivial relations between the three involutions $s_x,$ 
$s_y$ and $s_z,$ so that $\A$ is isomorphic to the free product of three copies of $\Z/2\Z.$ 

\vv

Summing up,  {\sl{the image of the mapping class group in $\Aut(S_{(A,B,C,D)}),$ the group
$\Gamma_2^*, $  and $\A$  correspond to finite index subgroups of the full automorphism 
group of $\Aut(S_{(A,B,C,D)}).$}}  This is the reason why we shall focus on the dynamics 
of $\Gamma_2^*=\A$ on the surfaces $S\in \SS.$


\subsection{Notations and remarks} 
The conjugacy class of a representation $\rho:\pi_1(\Sphere^2_4)\to \SL(2,\C)$ 
will be denoted $[\rho].$ In general, this conjugacy class is uniquely determined 
by its image $\chi(\rho)$ in  the character variety $\chi(\Sphere^2_4),$ and we shall identify $\chi(\rho)$ to $[\rho].$

Automorphisms of surfaces $S_{(A,B,C,D)}$ will be denoted by standard letters, like $f,$ $g,$ $h,$ ... ; the group $\A$ will be identified to its various realizations as subgroups
of $\Aut(S_{(A,B,C,D)}),$ where $(A,B,C,D)$ describes $\C^4.$
If $M$ is an element of $\Gamma_2^*,$ the automorphism 
associated to $M$ is denoted $f_M$; this provides an isomorphism between 
$\Gamma_2^*$ and each realization of $\A.$
If $f$ is an automorphism of $S_{(A,B,C,D})$ which is contained in $\A,$ 
$M_f$ will denote the unique element of $\Gamma_2^*$ which corresponds to $f.$ 

If $\Phi\in \MCG^*(\Sphere^2_4)$ is a mapping class, the associated 
automorphism of the character variety will be denoted by $f_\Phi.$ 

\vv

The character surfaces $S_D$ that appeared in the introduction in the
case of  the once punctured torus are
isomorphic to $S_{(0,0,0,D)}$ by a simultaneous multiplication of the variables
by $-1.$ As a consequence, the study of the dynamics on all 
character surfaces $S\in \SS$ contains the case of the once punctured torus.


\subsection{Dynamics at infinity}\label{par:prelim-di}

The group $\A$ also acts by birational transformations
of the compactification $\overline{S}$ of $S$ in $\P^3(\C).$ 
In this section, we describe the dynamics at infinity, i.e. on the
triangle $\Delta.$ 

If $f$ is an element of $\A,$ the birational transformation of ${\overline{S}}$
defined by $f$ is not everywhere defined. The set of its indeterminacy points 
is denoted by $\Ind(f)$; $f$ is said to be {\sl{algebraically stable}} if $f^n$ does
not contract any curve onto $\Ind(f)$ for $n\geq 0$ (see 
\cite{Sibony:Survey,Diller-Favre:2001} for this notion).

The group $\Gamma_2^*$ acts by isometries on  the Poincar\'e half plane $\H.$ 
Let $j_x,$ $j_y$ and $j_z$ be the three points on the boundary 
of $\H$ with coordinates $-1,$ $0,$ and $\infty$ respectively. 
The three generators  $r_x,$ $r_y,$ and $r_z$ of $\Gamma_2^*$ 
(see \ref{par:prelim-ca})  are  the reflections of  $\H$ around the three 
geodesics which join respectively  $j_y$ to $j_z,$  $j_z$ to $j_x,$ and $j_x$  
to $j_y$. As a consequence, $\Gamma_2^*$ coincides with the group of  symmetries
of the tesselation of $\H$ by ideal triangles, one of which has vertices $j_x,$ $j_y$ and $j_z.$ 
This picture will be useful to describe the action of $\A$ on $\Delta.$

\vv

First, one shows easily that the involution $s_x$  acts on the triangle $\Delta$ 
in the following way:  The image of the side $D_x$ is the vertex $v_x$ and the vertex $v_x$ is 
blown up onto the side $D_x$ ;  the sides $D_y$ and $D_z$ are invariant and $s_x$ 
permutes the vertices and fixes the middle points $m_y$ and $m_z$ of each of these sides. 
An analogous statement holds of course for $s_y$ and $s_z.$ In particular,  the
action of $\A$ at infinity does not depend on the  set of parameters $(A,B,C,D).$ 

Beside the three involutions $s_x,$ $s_y$ and $s_z,$ three new elements 
of $\A$ play a particular role in the study of the dynamics of $\MCG^*(\Sphere^2_4).$ 
These elements are
$$
g_x=s_z\circ s_y, \quad g_y=s_x\circ s_z , \quad  {\text{and}} \quad   g_z=s_y\circ s_x.
$$
They correspond to Dehn twists in the mapping class group. Each of them preserves 
one of the coordinate variables $x,$ $y$ or $z$ respectively. The
action of $g_x$ (resp. $g_y,$ resp. $g_z$) on $\Delta$ is the following: $g_x$ contracts
both $D_y$ and $D_z\setminus\{v_y \}$ on $v_z, $ and preserves $D_x$;  its inverse 
contracts $D_y$ and $D_z\setminus\{v_z\}$ on $v_y.$ In particular $\Ind(g_x)=v_y$ and 
$\Ind(g_x^{-1})=v_z.$ 

\vv

Let $f$ be any element of $\A\setminus \{{\text{Id}}\}$ and $M_f$ be the corresponding 
element of $\Gamma^*_2.$ If $M_f$ is elliptic, $f$ is conjugate to $s_x,$ $s_y$
or $s_z.$ If $M_f$ is parabolic, $f$ is conjugate to an iterate of $g_x,$ $g_y$ or
$g_z$ (see \cite{Cantat-Loray:2007}). In both cases, the action of $f$ on $\Delta$ 
has just been described. 

If $M_f$ is hyperbolic, the isometry $M_f$ of $\H$ has two fixed points 
at infinity, an attracting fixed point $\omega(f)$ and a repulsive fixed point $\alpha(f),$
and   the action of $f$ on $\Delta$ can be described as follows: The three sides
of $\Delta$ are blown down on the vertex $v_x$ (resp. $v_y$ resp. $v_z$) if
$\omega(f)$ is contained in the interval $[j_y,j_z]$ (resp. $[j_z,j_x],$ resp. $[j_x,j_y]$); 
the unique  indeterminacy point of $f$ is $v_x$ (resp. $v_y$ resp. $v_z$) if 
$\alpha(f)$ is contained in $[j_y, j_z]$ (resp. $[j_z,j_x],$ resp. $[j_x,j_y]$). In particular
$ \Ind(f)$ coincides with $\Ind(f^{-1})$ if and only if $\alpha(f)$ and $\omega(f)$
are in the same connected component of $\partial\H \setminus\{j_x, j_y, j_z\}.$
As a consequence, we get the following result (see \cite{Cantat-Loray:2007} for details). 

\begin{pro}\label{pro:IU}
Let $S$ be any member of the family $\SS.$ Let $f$ be an  element of $\A.$ Assume
that the element $M_f$ of $\Gamma_2^*$ that corresponds to $f$ is hyperbolic.
\begin{itemize}
\item The birational transformation $f:{\overline{S}}\to {\overline{S}}$ is algebraically 
stable if, and only if  $f$  is a cyclically reduced composition of the three involutions 
$s_x,$ $s_y$ and $s_z$ (in which each involution appears at least once). In particular,
any hyperbolic element $f$ of $\A$ is conjugate to an algebraically stable element of
$\A.$ 
\item If $f$ is algebraically stable, $f^n$ contracts the whole triangle 
$\Delta\setminus \Ind(f)$ onto $\Ind(f^{-1})$ as soon as $n$ is a  positive integer.
\end{itemize}
\end{pro}


\subsection{Topological entropy and types of automorphisms}

An element $f$ of $\A$ will be termed elliptic, parabolic or hyperbolic, according to the
type of the isometry $M_f\in \Gamma_2^*.$  By theorem B of \cite{Cantat-Loray:2007}, 
the {\sl{topological entropy}} $h_{top}(f)$ of  $f:S_{(A,B,C,D)}(\C)\to S_{(A,B,C,D)}(\C)$  
does not depend on the parameters $(A,B,C,D)$ and is equal to the logarithm of 
the spectral radius $\lambda(f)$ of $M_f$:
\begin{equation}
h_{top(f)}=\log (\lambda(f)).
\end{equation}
In particular, pseudo-Anosov mapping classes are exactly those with positive
entropy on the character surfaces $S_{(A,B,C,D)}.$ As explained in the previous
section, Dehn twists correspond to parabolic elements and are conjugate to a
power of $g_x,$ $g_y$ or $g_z,$ while elliptic automorphisms are conjugate
to $s_x,$ $s_y$ or $s_z.$ 

\begin{rem} This should be compared to the description of the group of
polynomial automorphisms of the affine plane $\C^2.$ If $h$ is an element
$\Aut(\C^2),$ either $h$ is conjugate to an elementary automorphism, which 
means that $h$ preserves the pencil of lines $y= c^{ste},$ or the topological 
entropy is equal to $\log (d(h)),$ where $d(h)$ is an integer (see \S 
\ref{par:boundedorbits-julia} for references).
\end{rem}


\subsection{The Cayley cubic}

The surface $S_4$ will play a central role in this paper. This surface is the 
unique element of $\SS$ with four singularities, and is therefore the unique 
element of $\SS$ that is isomorphic to the Cayley 
cubic (see \cite{Cantat-Loray:2007}). We shall call it "the Cayley cubic" 
and denote it by $S_C.$ This surface is isomorphic to the quotient of
$\C^*\times \C^*$ by the involution 
$\eta(x,y)=(x^{-1},y^{-1}).$ The map 
$$
\pi_C(u,v) = \left(u+\frac{1}{u},\, v+\frac{1}{v},\, uv + \frac{1}{uv}\right)
$$
gives an explicit isomorphism between $(\C^*\times \C^*)/\eta$ and $S_C$: Fixed points 
of $\eta,$ as $(-1,1),$ correspond to singular points of $S_C$.

The group $\GL(2,\Z)$ acts on $\C^*\times \C^*$ by {\sl{monomial transformations}}: 
If $M=(m_{ij})$ is an element of $\GL(2,\Z),$ and if $(u,v)$ is a point of $\C^*\times \C^*,$ 
then
$$
(u,v)^M= (u^{m_{11}}v^{m_{12}}, u^{m_{21}}v^{m_{22}}).
$$
This action commutes with $\eta$, so that $\PGL(2,\Z)$ acts on the quotient $S_C.$ 

The surface $S_C$ is one of the character surfaces for the once punctured torus: 
It corresponds to reducible representations of $\pi_1(\T_1)$ (with $\tr(\rho[\alpha,\beta])=2$). 
Of course, the monomial action of $\PGL(2,\Z)$ on $S_C$  coincides with the action
of the mapping class group of $\T_1$ on the character surface $S_C.$ Changing signs of
coordinates, it also coincides with the action of $\Gamma_2^*\subset \MCG(\Sphere^2_4)$ 
on the character surface corresponding to parameters $(a,b,c,d)=(0,0,0,0)$ or $(2,2,2,-2),$ 
up to permutation of $a,$ $b,$ $c,$ and $d$ and multiplication by $-1$ 
(see \cite{Cantat-Loray:2007}). 

The product $\C^*\times \C^*$ retracts by deformation onto the $2$-dimensional 
real torus $\Sphere^1\times \Sphere^1.$ The monomial action of $\GL(2,\Z)$ preserves
this torus: It acts "linearly" on this torus if we use the parametrization
$u=e^{2i\pi s},$ $v=e^{2i\pi t}.$  After deleting the four singularities of $S_C,$ the real 
part $S_C(\R)$ has five components, 
and the closure of the unique bounded component is 
the image of $\Sphere^1\times \Sphere^1$ by $\pi_C.$ The closure of the four unbounded components
are images of $\R^+\times\R^+,$ $\R^+\times\R^-,$ $\R^-\times\R^+,$ and $\R^-\times\R^-,$ 
by $\pi_C.$ 



\section{Elements with positive entropy}\label{par:Henon}

In this section, we describe the dynamics of hyperbolic elements in the group $\A$
on any surface $S_{(A,B,C,D)}(\C)$ of our family $\SS.$ 

Let $f$ be a hyperbolic element of $\A.$ 
After conjugation by an element $h$ of $\A,$ we can assume that $f$ is algebraically 
stable; in our context, this property means that, for any element $S$ of $\SS,$ the 
indeterminacy set of the birational transformation 
${\overline{f}}:{\overline{S}}\dasharrow {\overline{S}}$ and the indeterminacy set
of ${\overline{f}}^{-1}$ are two distinct 
vertices of the triangle at infinity $\Delta$ (see \S \ref{par:prelim-di}). 
In what follows, we shall assume that $f$ is  algebraically stable and denote
$\Ind(f^{-1})$ by $v_+$ and $\Ind(f)$ by $v_-.$


\subsection{Attracting basin of $\Ind(f^{-1})$}\label{par:attracting_bassin}


The birational transformation ${\overline{f}}$ is holomorphic in a neighborhood 
of $v_+$ and contracts $\Delta\setminus\{ v_-\}$ on $v_+.$ In particular, ${\overline{f}}$
contracts the two sides of $\Delta$ that contain $v_+$ on the vertex $v_+.$ 
Using the terminology of \cite{Favre:2000},  
${\overline{f}}$ determines a rigid and irreducible contracting germ near
$v_+.$ 

\vv

\begin{thm}
If $f$ is an algebraically stable hyperbolic element of $\A,$ then 
there exist an element $N_f$ of $\GL(2,\Z)$ with positive entries
which is conjugate to $M_f$ in $\PGL(2,\Z),$ a neighborhood ${\mathcal{U}}$ of $v_+$ 
in ${\overline{S}},$
and a  holomorphic diffeomorphism $\Psi^+_f:\disk \times \disk \to {\mathcal{U}}$ such
that $\Psi^+_f(0,0)=v_+$ and
$$
\Psi^+_f((u,v)^{N_f})=f(\Psi^+_f(u,v))
$$
for all $(u,v)$ in the bidisk $\disk \times \disk.$
\end{thm}

\begin{proof}
Let ${\mathcal{U}}$ be a small bidisk around $v_+,$ in which the two sides of $\Delta$
correspond to the two coordinate axis.  The fundamental group of ${\mathcal{U}}\setminus \Delta$ is isomorphic to $(\Z^2, +)$ and ${\overline{f}}$ induces an automorphism $N_f$ of this group. Since $f$ is a rigid and irreducible contracting germ near $v_+,$ a theorem of Dloussky and Favre asserts that $f$ is locally conjugate
to the monomial transformation that $N_f$ determines ; in particular, $f$ being a local contraction,  $N_f$ has
positive entries
(see class $6$ of the classification, Table II, and page 483 in \cite{Favre:2000}). 
The fact that the conjugacy  $\Psi_f$
is defined on the whole bidisk will be part of the next proposition.

To prove that $N_f$ is conjugate to $\pm M_f$ in $\GL(2,\Z),$ one argues as follows. 
The matrix $N_f$ is obtained from the action of $f$ on  the fundamental 
group of ${\mathcal{U}}\setminus \Delta.$ In the case of the Cayley cubic, 
$$
\pi_C: \C^*\times \C^* \to {\overline{S_C}}\setminus\Delta
$$ 
is a $2$ to $1$  covering,  $\C^*\times \C^*$ retracts by deformation on the
torus $\Sphere^1\times \Sphere^1,$ and the action of $f$ on the fundamental
group of  ${\mathcal{U}}\setminus \Delta$ is therefore covered by the action 
of $M_f$ on 
$
\pi_1( \Sphere^1\times \Sphere^1)=\Z\times \Z.
$
This implies that $N_f$ is conjugate to $M_f$ in $\PGL(2,Z).$ Since the general case is obtained
from the Cayley case by a smooth deformation, this is true for any set of
parameters $(A,B,C,D).$
\end{proof}

Let $s(f)$ be the slope of the eigenline of the linear planar transformation $N_f,$ 
which corresponds to the eigenvalue $1/\lambda(f)$;  $s(f)$ is a negative 
real number. The basin of attraction of the origin for the monomial transformation 
$N_f$ is
$$
{\overline{\Omega(N_f)}}=\{ (u,v)\in \C^2 \, \vert \quad  \vert v\vert < \vert u \vert^{s(f)}\}.
$$
In particular, this basin contains the full bidisk. We shall denote by $\Omega(N_f)$ the
intersection of ${\overline{\Omega(N_f)}}$ with $\C^*\times \C^*.$ 
Similar notations will be used for the basin of attraction   ${\overline{\Omega(v_+)}}$ 
of the point $v_+$ 
for ${\overline{f}}$ in ${\overline{S}},$ and for its intersection $\Omega(v_+)$  with $S.$ 

\begin{pro}\label{pro:conjugacy}
The conjugacy $\Psi^+_f$ extends to a  biholomorphism between  
$\Omega(N_f) $ and $\Omega(v_+).$ 
\end{pro}

\begin{proof}
Since the monomial transformation $N_f$ is contracting and $f:S\to S$ is invertible, 
we can extend $\Psi^+_f$ to $\Omega(N_f)\cap (\C^*\times \C^*)$ by the functional 
equation 
$$
\Psi^+_f(u,v)=f^{-n}(\Psi^+_f((u,v)^{N_f^{n}}),
$$
where $n$ is large enough for $(u,v)^{ N_f^{n} }$ to be in the initial domain of
definition of $\Psi^+_f.$ The map $\Psi^+_f:\Omega(N_f)\cap (\C^*\times \C^*)\to S$ is 
a local diffeomorphism, the image of which coincides with the basin of attraction of $v_+$ 
in $S.$ 
It remains to prove that the map $\Psi^+_f$ is injective. Assume  that $\Psi^+_f(u_1,v_1)=\Psi^+_f(u_2,v_2).$ Then
$f^n(\Psi^+_f(u_1,v_1))=f^n(\Psi^+_f(u_2,v_2)),$ and therefore 
$$
\Psi^+_f((u_1,v_1)^{N_f^n})=\Psi^+_f((u_1,v_1)^{N_f^n}),
$$
for any $n.$ Since $\Psi^+_f$ is injective in a neighborhood of the origin, and since the monomial
transformation $N_f$ is also injective, this implies $(u_1,v_1)=(u_2,v_2).$  \end{proof}

In what follows, $\Vert \, . \, \Vert$ will denote the usual euclidean norm in $\C^3.$

\begin{cor}
Let $f$ be an algebraically stable hyperbolic element of $\A.$ 
If $m$ is a point of $S$ with an unbounded positive orbit, then $f^n(m)$
goes to $\Ind(f^{-1})$ when $n$ goes to $+\infty$ and 
$$
\log\Vert f^n(m)\Vert \sim \lambda(f)^n.
$$ 
\end{cor}

\begin{proof} 
First we apply the previous results to the study of $f^{-1}$ and its basin 
of attraction near $v_{-}.$ Let us fix a small ball $B$ around $v_-$ in the surface 
${\overline{S}}.$ If $B$ is small enough, then $B$ is contained in the basin of attraction 
of $f^{-1}$: The orbit of a point $m_0\in B$ by $f^{-1}$ stays in $B$ and converges towards $v_-.$  Since $f$ contracts $\Delta\setminus\{v_-\}$ on $v_+,$ there is a 
neighborhood ${\mathcal{V}}\subset {\overline{S}} $ of $\Delta\setminus B$ which is contained in the basin of attraction of $v_+.$ 
Let $m$ be a point with unbounded orbit. Since ${\mathcal{V}}\cup B$ is a neighborhood of $\Delta,$
the sequence $(f^n(m))$ will visit ${\mathcal{V}}\cup B$ infinitely many times. Let $n_1$ be the first
positive time for which $f^{n_1}(m)$ is contained in ${\mathcal{V}}\cup B.$ 
Let $n_2$ be the first time after $n_1$ such that $f^{n_2}(m)$ escapes $B.$  
Then $f^n(m)$ never comes back in $B$ for $n>n_2.$ Pick a $n>n_2$ such that $
f^n(m)$ is contained in $ {\mathcal{V}}\cup B.$ Then $f^n(m)$ is in ${\mathcal{V}}$
and therefore in the basin of $v_+.$ This implies that the sequence 
$f^n(m)$ converges towards $v_+.$ 
In order to study the growth of $\Vert f^n(m)\Vert$ in a neighborhood of $v_+,$ 
 we apply the conjugacy
 $\Psi^+_f$: What we now need to control is the growth of $\Vert (u,v)^{N_f^n}\Vert^{-1},$ 
 and the result is  an easy exercise using exponential coordinates $(u,v)=(e^s, e^t),$ 
 in $\disc^*\times \disc^*.$  \end{proof}

\begin{cor}\label{cor:inv-curve}
If $f$ is a hyperbolic element of $\A$ and $A,$ $B,$ $C,$ and $D$ are four
complex numbers, $f$ does not preserve any algebraic
curve  in $S_{(A,B,C,D)}.$ 
\end{cor}

\begin{proof}
Let us assume the existence of a set of parameters $(A,B,C,D)$ and of an 
$f$-invariant algebraic curve $E\subset S_{(A,B,C,D)}.$ Let ${\overline{E}}$
be the Zariski-closure of $E$ in ${\overline{S}}_{(A,B,C,D)}(\C)$; $f$ induces an automorphism
${\overline{f}}$ of the compact Riemann surface ${\overline{E}}.$
Since $\C^3$ does not contain any $1$-dimensional compact complex 
subvariety, ${\overline{E}}$ contains points at infinity. These points must
coincide with $v_+$ and/or $v_-.$ In particular, the restriction of $f$ to ${\overline{E}}$
has at least one superattracting (or superrepulsive) fixed point. This is a contradiction 
with the fact that ${\overline{f}}:{\overline{E}}\to {\overline{E}}$ is an automorphism. 
\end{proof}


\subsection{Bounded orbits and Julia sets}\label{par:boundedorbits-julia}


Let us consider the case of a polynomial diffeomorphism $h$ of the affine plane $\C^2$
with positive topological entropy (an automorphism of  H\'enon type). 
After  conjugation by an element of $\Aut[\C^2],$ we may assume that $h$ is algebraically 
stable in $\P^2(\C).$ In that case, the dynamics of $h$ at infinity 
also exhibits two attracting fixed points, one  for $h,$ $w_+,$ and one for $h^{-1},$ 
$w_-,$ but there are three differences with the dynamics
of hyperbolic elements of $\A$: The exponential escape growth rate is an integer $d(h)$
(while $\lambda(f)$ is an irrational quadratic integer), the model 
to which $h$ is conjugate near $w_+$ is not invertible, 
and the conjugacy $\Psi_h$ is a covering map of infinite degree between the basins of attraction. 
We refer the reader to \cite{HOV:1994}, \cite{Favre:2000} and
\cite{HPV:2000} for an extensive study of this situation.

Beside these differences, we shall see that the dynamics of hyperbolic elements of $\A$ 
is similar to the dynamics of H\'enon automorphisms. 
In analogy with the H\'enon case, let us introduce the following definitions:

\vspace{0.1cm}

 $\bullet$ $K^+(f)$ is the set of bounded forward orbits. This is also the set 
of points $m$ in the surface $S,$ for which $(f^n(m))$ does not converge
to $v_+$  when $n$ goes to $+\infty.$ 

$K^-(f)$ is the set of bounded backward orbits, and 
$K(f)=K^-(f)\cap K^-(f).$

\vspace{0.1cm}

 $\bullet$ $J^+(f)$ is the boundary of $K^+(f),$ $J^-(f)$ is the boundary of $K^-(f),$
and  $J(f)$ is the subset of $\partial K(f)$ defined by
$
J(f)=J^-(f)\cap J^+(f).
$
The set $J(f)$ will be called {\sl{the Julia set of $f$}}.

\vspace{0.1cm}

$\bullet$  $J^*(f)$ is the closure of the set of saddle periodic points of $f$ (see below).

\vspace{0.1cm}

Figure 1, right,  is  a one dimensional (complex) slice of $K^+(f)$ when 
$f= s_x\circ s_y\circ s_z$ and
$
(A,B,C,D)=(0,0,0,0).
$
 The left part of this figure 
represents a few hundred orbits of this automorphism on the unique 
compact connected component  of $S_{(0,0,0,2)}(\R).$ 


\subsection{Green functions and dynamics}\label{par:henonlike} 


We define the Green functions of $f$ by 
\begin{eqnarray}\label{eq:green}
G_f^+(m) & = &\lim_{n\to +\infty} \frac{1}{\lambda(f)^n}\log^+ \Vert f^n(m)\Vert, \\
G_f^-(m) & = & \lim_{n\to +\infty} \frac{1}{\lambda(f)^n}\log^+ \Vert f^{-n}(m)\Vert.
\end{eqnarray}
By proposition \ref{pro:conjugacy} and its corollary, both functions are well defined
and the zero set of $G_f^\pm$ coincides with $K^\pm(f).$ Moreover, the convergence
is uniform on compact subsets of $S.$ Since $\log^+ \Vert \, \Vert$
is a pluri-subharmonic function, $G_f^+$ (resp. $G_f^-$) is pluri-subharmonic and is
pluri-harmonic on the complement of $K^+(f)$ (resp. $K^-(f)$) (see \cite{BS:1991inv, Fornaess-Sibony:1992, Sibony:Survey}
for the details of the proof). These functions satisfy the invariance properties 
\begin{equation}
G^+_f\circ f = \lambda(f) G^+_f \quad {\text{and}} \quad 
G^-_f\circ f = \lambda(f)^{-1} G^-_f
\end{equation}
The following results have been proved for H\'enon mappings;  we list 
them with appropriate references, in which the reader can find a proof which
applies to our context (see also \cite{Cantat:Acta}, \cite{Bedford-Diller:2005}, 
\cite{Dujardin:2006},  \cite{Sibony:Survey} for similar contexts). 

\vspace{0.1cm}

$\bullet$ $G_f^+$ and $G_f^-$ are H\"older continuous (see \cite{Dinh-Sibony:AMS}, sections 2.2, 2.3). 
The currents 
\begin{equation} 
T_f^+=dd^cG_f^+ \quad {\text{and}} \quad T_f^-=dd^cG_f^-
\end{equation}
are closed and positive, and $f^*T^\pm_f=\lambda(f)^\pm T^\pm_f.$ 
By \cite{BS:1991inv}, section 3, the support of $T^+_f$ is $J^+(f),$ the support of $T^-_f$ is $J^-(f)$  (see also \cite{Sibony:Survey}).

\vspace{0.1cm}

$\bullet$ Since the potentials $G^+_f$ and $G^-_f$ are continuous, the product 
\begin{equation} 
\mu_f=T^+\wedge T^-
\end{equation} 
is a well defined positive measure, and is $f$-invariant.
Multiplying $G^+_f$ and $G^-_f$ by positive constants, we can, and we shall assume that 
$\mu_f$ is a probability measure. \hfill (see \cite{BS:1991inv},~section 3)

\vspace{0.1cm}

$\bullet$ The topological entropy of $f$ is $\log(\lambda(f))$ and the measure $\mu_f$ is the unique $f$-invariant probability measure with maximal entropy. (see \cite{BLS:I},~section 3)

\vspace{0.1cm}

$\bullet$ If $m$ is a saddle periodic point of $f,$ its unstable (resp. stable)
manifold $W^u(m)$ (resp. $W^s(m)$) is parametrized by $\C.$ Let $\xi:\C\to S$
be such a parametrization of $W^u(m)$ with $\xi(0)=m.$ Let $\disk\subset \C$ be the unit
disk, and let $\chi$ be a smooth non negative function on $\xi(\disk),$ with 
$\chi(m)>0$ and $\chi=0$   in a neighborhood of $\xi(\partial\disk).$ Let $[\xi(\disk)]$ be 
the current of integration on $\xi(\disk).$ The sequence of currents
$$
\frac{1}{\lambda(f)^n}f^{-n}_*(\chi . [\xi(\disk)])
$$
weakly converges toward a positive multiple of $T^-_f.$ (see \cite{BS:1991ams},
sections 2 and  3, \cite{Fornaess-Sibony:1992})

\vspace{0.1cm}

$\bullet$ By corollary \ref{cor:inv-curve},  periodic points of $f$ are isolated.  
The number of periodic points of period $N$ grows like $\lambda(f)^N.$  Most of them  are hyperbolic saddle points: 
If ${\mathcal{P}}(f,N)$ denotes either the
set of  periodic points with period $N$ or the set of 
periodic saddle points of period $N,$ then 
$$
\sum_{m\in {\mathcal{P}}(f,N)} \delta_m \to \mu_f
$$
where the convergence is a weak convergence in the space of probability 
measures on compact subsets of $S.$ (see \cite{BLS:I}, \cite{BLS:II}, 
and \cite{Dujardin:2006})

\vspace{0.1cm}

$\bullet$ The support $J^*(f)$ of $\mu_f$ simultaneously coincides 
with the Shilov boundary of $K(f)$ and with the closure of 
periodic saddle points of $f.$ In particular, any periodic saddle 
point of $f$ is in the support of $\mu_f.$ If $p$ and $q$ are periodic
saddle points, then $J^*(f)$ coincides with the closure of
$W^u(p)\cap W^s(q).$ (see \cite{BLS:I} and \cite{BLS:II})

\vspace{0.1cm}

$\bullet$ Since $f$ is area preserving (see \S \ref{par:prelim-amg}), the interior 
of $K(f),$ $K^+(f)$ and $K^-(f)$ coincide. In particular, the interior of $K^+(f)$
is a bounded open subset of $S(\C).$ (see lemma 5.5 of \cite{BS:1991inv}) 


\section{The quasi-fuchsian locus and its complement}\label{par:quasifuchsian}


In this section, we shall mostly restrict the study to the case of the once punctured
torus with a cusp, and provide hints for more general statements. 


\subsection{Quasi-fuchian space and Bers' parametrization}

Let $\T_1$ be a once punctured torus.  Let  $\Te(\T_1)$ be the
 Teichm\"uller space of  complete hyperbolic metrics on $\T_1$ with finite
area $2\pi,$ or equivalently with a cusp at the puncture: $\Te(\T_1)$ is isomorphic, 
and will be identified, to the upper half plane $\H^+.$ The dynamics of 
$\MCG(\T_1)$ on $\Te(\T_1)$ is conjugate to the usual  action of 
$\PSL(2,\Z)$ on $\H^+.$

Any point in the Teichm\"uller space gives rise to a representation
${\overline{\rho}}:F_2\to \PSL(2,\R)$ that can be lifted to four distinct representations
into $\SL(2,\R).$ The cusp condition gives rise to the same equation 
$\tr(\rho[\alpha,\beta])=-2$
for any of these four representations. This provides four embeddings of the
Teichm\"uller space into the surface $S_0(\R)$: The four images are the 
four unbounded  components of $S_0(\R),$ each of which is diffeomorphic to 
$\H^+$; apart from these four components, $S_0(\R)$ 
contains an isolated singularity at the origin. This point corresponds to the 
conjugacy class of the representation $\rho_{quat},$ defined by
\begin{equation}\label{eq:quat}
\rho_{quat}(\alpha)=\left( 
\begin{array}{cc} 0 & i \\ i & 0 \end{array}\right), 
\quad \rho_{quat}(\beta)=
\left( 
\begin{array}{cc} 0 & -1 \\ 1 & 0 \end{array}\right).
\end{equation}
Its image coincides with the quaternionic group of order eight.
The mapping class group of the torus 
acts on $S_0(\R),$ preserves the origin and the connected component 
$$
S_0^+(\R)=S_0(\R)\cap (\R^+)^3,
$$
and permutes the remaining three components.

Let $\DF\subset S_0(\C) $ be the set of conjugacy classes of discrete and faithful 
representations $\rho:F_2\to \SL(2,\C)$ with $\tr(\rho[\alpha,\beta])=-2.$  
This set is composed of four distinct connected components, one of
them, $\DF^+,$ containing $S_0^+(\R).$
The component $S_0^+(\R)$ is made of conjugacy classes of fuchsian 
representations, and the set $\QF$ of their quasi-fuchsian deformations coincides with 
the interior of $\DF^+$ (see \cite{Minsky:Survey}, and the references therein). 

Let $\T_1'$ be the once punctured torus with the opposite orientation. 
Bers' parametrization of the space of quasi-fuchsian 
representations provides a holomorphic bijection 
$$
\Be : \Te(\T_1)\times  \Te(\T_1') \to Int(\DF^+).
$$
We may identify $\Te(\T_1)$ with the upper half plane $\H^+$ and
$\Te(\T_1')$ with the lower half plane $\H^-.$ The group $\PSL(2,\Z)$
acts on $\P^1(\C),$ preserving $\P^1(\R),$ $\H^+,$ and $\H^-.$ In 
particular, $\MCG(\T_1)=\SL(2,\Z)$ acts diagonally on 
$$
\Te(\T_1)\times \Te(\T_1')=\H^+\times \H^-.
$$
With these identifications, the map $\Be$ conjugates the diagonal action of 
$\MCG(\T_1)$ on $\H^+\times \H^-$ with its action on the character variety: If $\Phi$
is a mapping class and $f_\Phi$ is the automorphism of $S_0$ which is
determined by $\Phi,$ then 
$$
\Be(\Phi(X),\Phi(Y)) = f_\Phi(\Be(X,Y))
$$
for any $(X,Y)$ in $\H^+\times \H^-.$ It conjugates the action of 
$\MCG(\T_1)$ on the set  
$$
\{(z_1,z_2)\in \H^+\times \H^-\vert z_1={\overline{z_2}}\}
$$
with the corresponding action on $S_0^+(\R).$ 
The Bers map  extends up to the boundary of $\H^+\times \H^-$ minus
its diagonal (we shall call it  the restricted boundary, and denote it
by $\partial^*(\H^+\times \H^-)$).
Minsky proved in \cite{Minsky:1999} that $\Be$ induces a continuous 
bijection from $\partial^*(\H^+\times \H^-)$ to 
the boundary of $\DF^+.$


\subsection{Mapping torus and fixed points}\label{par:qffp}

Let $\Phi\in \MCG(\T_1)$ be a pseudo-Anosov mapping class. Let $X_\Phi$ be the 
mapping torus determined by $\Phi$:  The threefold $X_\Phi$ is  obtained 
by suspension of $\T_1$ over the circle, with monodromy $\Phi.$ 
Thurston's hyperbolization theorem tells us that $X_\Phi$ can be endowed with 
a complete hyperbolic metric of finite volume. This provides a discrete 
and faithful  representation 
$$
\rho_\Phi:\pi_1(X_f)\to {\sf{Isom}}(\H^3)=\PSL(2,\C)
$$
If we restrict $\rho_\Phi$ to the fundamental group  of  the torus fiber 
of $X_\Phi,$ and if we choose the appropriate lift to $\SL(2,\C),$ 
we get a point $[\rho_\Phi]$ in $\DF^+\subset S_0(\C)$ which is fixed by 
the automorphism $f_\Phi.$  
Let $\alpha(\Phi)$  (resp. $\omega(\Phi)$) be the repulsive (resp. attracting) fixed point 
of $\Phi$ on the boundary of $\Te(\T_1).$ Since $\Be$ is a continuous conjugacy, we have
$$
{\Be}(\alpha(\Phi),\omega(\Phi))= [\rho_\Phi].
$$
The fixed point $(\omega(\Phi),\alpha(\Phi))$ provides a second fixed point on the
boundary of $\DF^+$: This point may be obtained by the same construction 
with $\Phi^{-1}$ in place of $\Phi.$ 
In \cite{McMullen:book}, McMullen proved that $[\rho_\Phi]$ is  a hyperbolic fixed 
point of $f_\Phi.$ The stable and unstable manifolds of $f_\Phi$ at $[\rho_\Phi]$ 
intersect $\DF^+$ along its boundary, 
\begin{eqnarray}
W^u([\rho_\Phi])\cap \DF^+ & = & \Be (\{ \alpha(\Phi) \}\times {\overline{\H^-}}\setminus \{(\alpha(\Phi),\alpha(\Phi))\}), \\
W^s([\rho_\Phi])\cap \DF^+ & =  & \Be ({\overline{\H^+}}\times \{\omega(\Phi)\} \setminus \{(\omega(\Phi),\omega(\Phi))\}).
\end{eqnarray}
In particular, the union of all stable manifolds $W^s([\rho_\Phi])\cap \DF^+$ of $f_\Phi,$ where $\Phi$ describes the set of pseudo-Anosov mapping classes, form a dense subset
of $\partial \DF^+.$ 

\begin{rem}
Each pseudo-Anosov class $\Phi$ determines an automorphism $f_\Phi,$ and therefore
a subset $K^+(f_\Phi)$ of $S_0(\C).$ The complement $\Omega^+(f_\Phi)$ of 
$K^+(f_\Phi)$ is open: It co\"incides  with the bassin of attraction of $f_\Phi$ at infinity. The interior of $\DF^+$ is contained in the intersection 
$$
\Omega(\MCG(\T _1)) := \bigcap_\Phi \, \Omega^+(f_\Phi)
$$
where $\Phi$ describes the set of pseudo-Anosov classes in the mapping class group
$\MCG(\T_1)=\SL(2,\Z).$ Since stable manifolds are dense in the boundary of $\DF,$
the quasi fuchsian locus $\Int(\DF^+)$ is a connected component of the interior of
$\Omega(\MCG(\T _1)).$ 
\end{rem}

\subsection{Two examples}\label{par:qfex}

We now present  two orbits in the complement of $\DF.$

\begin{thm} Let $\Phi$ be any pseudo-Anosov mapping class and $[\rho_\Phi]$ 
be one of the two fixed points of $f_\Phi$ on the boundary of the subset $\DF^+$ of $S_0(\C).$ 
There exists a representation $\rho_0: \pi_1(\T_1)\to \SL(2,\C)$ such that
 $[\rho_0]$ is an element of $S_0(\C)$ and
\begin{itemize}
\item the sequence $(f_\Phi)^n[\rho_0]$ converges toward the discrete and 
faithful representation $[\rho_\Phi]$ when $n$ goes to $+\infty$;
\item the closure of the mapping-class group orbit of $[\rho_0]$ contains
the origin $(0,0,0)$ of $S_0(\C),$ i.e. the conjugacy class of the finite representation 
$\rho_{quat}.$ 
\end{itemize} 
\end{thm}

\begin{proof} The fixed point $[\rho_\Phi]$ is hyperbolic, with a stable
manifold $W^s([\rho_\Phi]).$ The origin $(0,0,0)$ is the unique singular 
point of $S_0(\C).$ It corresponds to the representation $[\rho_{quat}]$ which 
is defined by equation (\ref{eq:quat}). This point is fixed by $f,$ and a direct computation shows that the differential of $f$ at the origin has finite order (order $1$ or $2$).

From section \ref{par:henonlike}, the interior of $K^+(f)$ coincides with the interior of $K^-(f)$ and is therefore an $f$-invariant bounded open subset of $S_0(\C).$ In particular, $\Int(K^+(f))$ is Kobayashi hyperbolic, and if $[\rho_{q}]$ is in $\Int(K^+(f)),$ then 
$f$ is locally linearizable around the origin $[\rho_{q}].$ Since $Df_{[\rho_{q}]}$ has finite order, $f$ would have finite order too. This contradiction shows that $[\rho_{q}]$ is not 
in the interior of $K^+(f).$

According to a theorem  of Bowditch (see theorem 5.5 of \cite{Bowditch:1998}), there exists
a neighborhood $U$ of the origin in $S_0(\C)$ with the property that any 
mapping class group orbit starting in $U$ contains the origin in its closure.

We know that $W^s([\rho_f])$ is dense in 
the boundary of $K^+(f)$ (see \S \ref{par:henonlike}).
Since $[\rho_{q}]$ is in $\partial K^+(f),$ $W^s([\rho_f])$
intersects the Bowditch's neighborhood $U.$ Any point $[\rho_0]$ in $W^s([\rho_f]) \cap U$
satisfies the properties of the theorem. 
\end{proof}

\begin{proof}[Proof of theorem \ref{thm:dg}] Let us now consider the 
dynamics of the mapping class 
$$
\Psi = \left(\begin{array}{cc} 2 & 1 \\ 1 & 1\end{array}\right)
$$
on the surface $S_2(\C).$ The surface $S_2(\C)$  corresponds to representations
$\rho:G\to \SL(2,\C),$ where $G=\langle \alpha, \beta\vert [\alpha,\beta]^4\rangle$
(see \S \ref{par:intro-gold}). 

As explained, for example, in \cite{McMullen:book}, section 3.7, the surface $S_2(\C)$ contains
an $f_\Psi$-invariant open subset corresponding to quasi-fuchsian deformations of 
the fuchsian groups obtained by endowing a hyperbolic metric on $\T_1$ with an
orbifold point of angle $\pi$ at the puncture. Thurston's hyperbolization theorem
provides a hyperbolic fixed point $[\rho_\Psi]$  of
$f_\Psi$ on the boundary of this set: The representation $\rho_\Psi:G\to \SL(2,\C)$ 
is discrete and faithful and comes from the existence of a hyperbolic structure on the 
complement of the figure eight knot, with an orbifold structure along the knot. 

The subset of $S_2(\C)$ corresponding to conjugacy classes of $\SU(2)$
representations coincides with the unique bounded connected component of 
$S_2(\R),$ and is homeomorphic to a sphere (see \cite{Goldman:1988}, figure 4). This
component is $f_\Psi$-invariant, and $f_\Psi$ has exactly two fixed points on it, namely 
$(x,x/(x-1),x),$ with
$$
x= \frac{\sqrt{17} + \sqrt{1+\sqrt{17}/2}}{2} \quad {\text{ or }} \quad \frac{\sqrt{17} - \sqrt{1+\sqrt{17}/2}}{2}.
$$
Both of them are saddle points. Let $[\rho_{\SU}]$ be one of these fixed points, and
let $W^s([\rho_\SU])$ and $W^u([\rho_\Psi])$ be the stable and unstable manifolds of
$f_\Psi$ through   $[\rho_{\SU}].$ 

From section \ref{par:henonlike}, 
we know that  $W^s([\rho_\SU])$ intersects $W^u([\rho_\Psi]).$ Let $[\rho_0]$
be one of these intersection points. The $f_\Psi$-orbit of $[\rho_0]$
contains both $[\rho_\Psi]$ and $[\rho_\SU].$ 

Finite orbits of $\MCG(\T_1)$ are listed in \cite{Dubrovin-Mazzocco:2000} and
$[\rho_\SU]$ does not appear in the list. As a consequence, the mapping class group 
orbit of $[\rho_\SU]$ is infinite and dense in the component of $\SU(2)$-representations 
(see \cite{Goldman:1997}, \cite{Goldman:2003}, or \cite{Cantat-Loray:2007}). This 
implies that the closure of the orbit of $[\rho_\SU]$ contains both $[\rho_\Psi]$ and the
$\SU(2)$-component of $S_2(\R).$ 
\end{proof}


\section{Real Dynamics of Hyperbolic Elements}\label{par:realdynamics}

In this section, we study the dynamics of hyperbolic elements on the real 
surfaces $S_{(A,B,C,D)}(\R)$ when the parameters are real numbers. 
The main goal of this section is to prove theorem \ref{thm:hypreal+} below, 
which extends, and precises,  theorem \ref{thm:hypreal}.

\subsection{Maximal entropy}


Let us fix a hyperbolic element $f\in \A.$ If 
the parameters $(A,B,C,D)$ are real, we get two dynamical systems: The first one
takes place on the complex surface $S(\C)$ and its main stochastic
properties have been listed in section \ref{par:henonlike}; the second one is induced
by the restriction of $f$ to the real part $S(\R).$ From time to time, we shall use the 
notation $f_\R$ for the restriction of $f$ to $S(\R).$ For example, we shall say that 
{\sl{$f_\R$ has maximal entropy}} if the entropy of $f:S(\R)\to S(\R)$
is equal to the topological entropy of $f:S(\C)\to S(\C),$ {\sl{i.e.}} to $\log(\lambda(f)).$ 

\begin{thm}\label{thm:meco}
Let $f$ be a hyperbolic element of $\A.$ If $(A,B,C,D)$ are real parameters, the
following conditions are equivalent: 
\begin{enumerate}
 \item $f_\R$ has maximal entropy;
\item $J^*(f)$ is contained in $S(\R)$;
\item $K(f)$ is contained in $S(\R).$
\end{enumerate}
In that case, $J^*(f)=J(f)=K(f).$
\end{thm}

This theorem is an easy consequence of the results of  section \ref{par:henonlike}
(see \cite{BLS:I}, section 10 for a proof). Our first goal is to prove the following result.

\begin{thm}\label{thm:entropy_minoration}
Let $f$ be a hyperbolic element of $\A.$ If $(A,B,C,D)$ are real parameters
such that $S_{(A,B,C,D)}(\R)$ is connected, then $f_\R$ has maximal entropy.
\end{thm}

\begin{rem}\label{rem:bgco}
Benedetto and Goldman studied the various topologies that can occur for
$S(\R).$ 
Using $(a,b,c,d)$ parameters (see section \ref{par:prelim-equations}),
 $S(\R)$ is connected if and only  $(i)$ none of the parameters $a,$ $b,$ $c,$ and $d$
is contained in the interval $(-2,2)$ and $(ii)$ the product $abcd$ is negative. In that
case, the surface $S(\R)$ is homeomorphic to a sphere minus four punctures 
(see \cite{Benedetto-Goldman:1999}). These conditions on $(a,b,c,d)$ define four arcwise connected subsets of
$\R^4,$ that contain respectively the $8$ points 
$( 2\epsilon_1, 2\epsilon_2,  2\epsilon_3, 2 \epsilon_4),$ with $\epsilon_i=\pm 1$ and $\Pi \epsilon_i=-1.$ 
All  these points correspond to the same surface $S_{(0,0,0,4)}$, {\sl{i.e.}} to the Cayley cubic
$S_C.$ As a consequence, any connected surface $S(\R)$ can be smoothly deformed 
to the Cayley cubic $S_C$ inside $\SS.$
\end{rem}

Before giving a proof of  theorem \ref{thm:entropy_minoration}, let us review a result of Bowen
concerning topological lower bounds for the entropy (see \cite{Bowen:1978}). 
Let $f$ be a homeomorphism of a marked topological space $(X,m),$ by which 
we mean that $m$ is a fixed point of $f.$  
Then, $f$ determines an automorphism 
$f_*:\pi_1(X,m)\to \pi_1(X,m).$ 
Let us assume that $\pi_1(X)$ is finitely generated, and fix a finite set  $\{\alpha_1, ..., \alpha_k\}$ 
of generators for $\pi_1(X).$ The growth rate of $f_*$ is defined to be 
$$
\lambda(f_*)=\limsup_{n\to +\infty} \left(\frac{1}{n} \diam (f^n(B))\right)
$$
where $\diam$ is the diameter with respect to the word metric (using the generators
$\alpha_i$) and $B$ is the ball of radius $1$ with respect to this metric. 
Bowen's theorem shows that 
$$
h_{top}(f)\geq \log (\lambda(f_*))
$$
as soon as $f$ is a continuous transformation of a compact manifold. 
Even though   $S(\R)$ is not compact, we can apply this theorem because unbounded 
orbits are contained in the basins of attraction
of $\Ind(f^{-1})$ and $\Ind(f).$

\begin{proof}[Proof of theorem \ref{thm:entropy_minoration}]
Let us first study  the case of the Cayley cubic $S_C.$ This surface is singular, 
and $S_C(\R)\setminus \Sing(S_C)$ contains a unique bounded component. This component
$S_C(\R)^0$ is a sphere with four punctures and the dynamics of $\A$ ({\sl{i.e.}} $\Gamma_2^*$) 
is covered by the monomial action of $\Gamma_2^*$ on the torus $\Sphere^1\times \Sphere^1$
in $\C^*\times \C^*.$ As a consequence, for any hyperbolic element $f$ in $\Gamma_2^*,$
the entropy of $f$ on $S_C(\R)^0$ is maximal. Moreover, in that case, the expanding factor 
$\lambda(f_*)$ coincides with the dynamical degree $\lambda(f),$ and Bowen's 
inequality is an equality.

If we deform the Cayley cubic in such a way that the surface $S(\R)$ is smooth and
connected, then $S(\R)$ is homeomorphic to a four punctured sphere (the punctures are now at
infinity), and the action of $f$ on the fundamental group 
of $S(\R)$ has not been changed along the deformation. As a consequence, 
Bowen's inequality gives 
$$
h_{top}(f_\R)\geq \log (\lambda(f))
$$
and the conclusion follows from 
$$h_{top}(f_\R)\leq h_{top}(f_\C)= \log(\lambda(f)).
$$ 
This concludes the proof
 for  smooth and connected surfaces $S(\R)$ (see remark \ref{rem:bgco}).
If $S(\R)$ is not smooth but is connected, then $S(\R)$ is a limit of smooth
connected members of the family $\SS.$ By semicontinuity of topological entropy, 
$f_\R$ has maximal entropy (see \cite{Newhouse:1989}). 
\end{proof}

\begin{cor}
Let $a,$ $b,$ $c,$ and $d$ be four real parameters in $\R\setminus [-2,2],$ the product of which
is negative. Let  $\rho:\pi_1(\Sphere_4)\to \SL(2,C)$ be a representation with 
boundary traces $a,$ $b,$ $c,$ and $d.$  Let $\Phi\in \Aut(\pi_1(\Sphere^2_4))$ 
be a pseudo-Anosov automorphism.
If $\rho\circ \Phi$ is conjugate to $\rho,$ 
then $\rho$ is  conjugate to a representation into $\SL(2,\R).$
 \end{cor}
 
\begin{proof}
Let $S$ be the element of the family $\SS$ that corresponds to the parameters $(a,b,c,d).$
The assumption on the parameters $a,$ $b,$ $c,$ and $d$ implies
that  $S(\R)$ is connected 
(see remark \ref{rem:bgco}), and that there is no $\SU(2)$-component
(this is obvious if $S(\R)$ is smooth, since $\SU(2)$ representations
would form a compact component, and this follows from 
\cite{Benedetto-Goldman:1999} in the singular case). 

If $\rho\circ \Phi^{-1}$ is conjugate to $\rho,$ then $\chi(\rho)$ is a fixed point of the automorphism
$f_\Phi$ induced by $\Phi$ on the surface $S.$ Since $S(\R)$ is connected,  $f_\R$ has maximal entropy. 
By theorem \ref{thm:meco}, all periodic points of $f$ are contained in $S(\R).$ 
This implies that $\rho$ is conjugate to an $\SL(2,\R)$-valued representation. \end{proof}


\subsection{Maximal entropy and quasi-hyperbolicity}\label{par:rdmeqh}


Bedford and Smillie recently developped a nice theory for H\'enon transformations
which extends the notion of quasi-hyperbolicity, a notion that had been previously introduced
for the dynamics of rational maps of one  complex variable. 
This theory can be applied to our context in order to study hyperbolic
automorphisms with maximal entropy. 

\subsubsection{Quasi-hyperbolicity}
Let $Sadd(f)$ be either the set of periodic saddle points of $f$ or the
set $W^u(p)\cap W^s(q)$ where $p$ and $q$ are two periodic fixed points
of $f$ (see \cite{BS:2002} for more examples). With such a choice, $Sadd(f)$ is $f$-invariant 
and its closure coincides with $J^*(f)$ (see \S \ref{par:henonlike}). Each point $m$ of 
$Sadd(f)$ has a stable manifold $W^s(p)$ and an unstable manifold
$W^u(p),$ and we can find two injective immersions $\xi_m^u,$ $\xi_m^s:\C\to S$ 
such that $\xi_m^{u/s}(0)=m,$ $\xi_m^{u/s}(\C)=W^{u/s}(m),$ and 
$$
\max \{ G^{+/-}(\xi_m^{u/s}(t))\, \vert \quad t \in \disk\} =1,
$$
where $\disk$ is the unit disk.
The parametrization $\xi^u_m$ and $\xi^s_m$ are uniquely determined by this
normalization up to a rotation of $t.$ Since $Sadd(f)$ is $f$-invariant 
and $f$ sends the unstable manifold at $m$ on the 
unstable manifold at $f(m),$ there is a non zero complex
number $\lambda(m)$ such that 
$$
f(\xi^u_m(t))=\xi^u_{f(m)}(\lambda(m)t).
$$
The number $\lambda(m)$ depends on the choices made for $\xi^u_m$ and $\xi^u_{f(m)}$
but its modulus $\vert \lambda(m) \vert$ only depends on $m.$ Since $G^+\circ f = \lambda(f)G^+,$
we obtain easily the inequality 
$\vert\lambda(m)\vert >1$ for all $m \in Sadd(f).$

We shall also need the growth function $\gro_m(r)$ of $G^+$ along the unstable
manifold $W^u(m),$ which is defined by 
$
\gro_m(r)=\max_{\vert t \vert \leq r } \left\{ G^+(\xi^u_m(t))\right\} ,
$ 
and the uniform growth function 
$$
\Gro(r) = \sup_{m\in Sadd(f)} \left\{\gro_m(r)\right\}.
$$
Bedford and Smillie proved in  \cite{BS:2002}, section 1, that
the following properties are equivalent:
\begin{enumerate}
 \item the family $\{\xi^u_m\, \vert \, m\in Sadd(f)\}$ is a normal 
family;
\item $\Gro (r_0)< \infty$ for some $1< r_0 < \infty$;
\item there exists $\kappa>1$ such that $\vert \lambda(m)\vert \geq \kappa$ for all
$m$ in $Sadd(f)$;
\item $\exists \,C,\beta<\infty$ such that $\gro_m(r)\leq Cr^\beta$ for all $m$ in $S$
and $r\geq 1.$
\end{enumerate} 
If one, and then all, of these properties is satisfied, $f$ is said to be {\sl{quasi-expanding}}.
If $f$ and $f^{-1}$ are quasi-expanding, then $f$ is said to be {\sl{quasi-hyperbolic}}.

\subsubsection{Maximal entropy}
It turns out that real H\'enon mappings with maximal entropy are necessarily quasi-hyperbolic
(see \cite{BS:2002}, theorem 4.8 and proposition 4.9). The proof of this result can be applied 
word by word to our context, and gives rise to the following theorem. 

\begin{thm}[Bedford Smillie, \cite{BS:2002} and \cite{BS:2004}]\label{thm:bsrh}
Let $f$ be a hyperbolic element of $\A$ and $S$ be an element of $\SS$ defined
by real parameters $(A,B,C,D).$ If $f_\R$ has maximal entropy, then $f$ is
quasi-hyperbolic, and any periodic point $m$ of $f$ 
is a saddle point, with $\vert \lambda(m) \vert \geq \lambda(f).$
\end{thm}

\begin{cor}
Let $f$ be a hyperbolic element of $\A$ and $S$ be an element of $\SS$ defined
by real parameters $(A,B,C,D).$ If $S(\R)$ is connected, then $f_\R$ has
maximal entropy and is quasi-hyperbolic. 
\end{cor}

\subsubsection{Uniform hyperbolicity and consequences}
In a subsequent paper, Bedford and Smillie also obtain a precise obstruction 
to uniform hyperbolicity. Let $p\in S(\R)$ be a saddle periodic point of $f.$
The unstable manifold of $p$ in $S(\R)$ is the intersection of $S(\R)$ with
the complex unstable manifold $W^u(p).$ This real unstable manifold is diffeomorphic to the
real  line $\R,$ and $p$ disconnects it into two half lines. If one of these half unstable manifolds
is contained in the complement of $K^+(f),$ one says that $p$ is {\sl{$u$-one-sided}}; 
 {\sl{$s$-one-sided}} points are defined in a similar way.  

\begin{thm}[Bedford Smillie, \cite{BS:2004}]\label{thm:bscme}
Let $f$ be a hyperbolic element of $\A$ and $S$ be an element of $\SS$ defined
by real parameters $(A,B,C,D).$ If $f_\R$ has maximal entropy but $K(f)$
is not a hyperbolic set for $f,$ then
\begin{itemize}
 \item there are periodic saddle points $p$ and $q$ (not necessarily distinct)
so that $W^u(p)$  intersects $W^s(q)$ tangentially with order 2 contact~;
 \item $p$ is $s$-one-sided and $q$ is $u$-one-sided~;
\item the restriction of $f$ to $K(f)$ is not expansive. 
\end{itemize}\end{thm}

\begin{thm}
Let $f$ be a hyperbolic element of $\A.$ Let $S$ be a smooth surface in the family
$\SS$ which is defined by real parameters $(A,B,C,D).$ If one of the 
connected components of $S(\R)$ is bounded, then the entropy of $f_\R$ is
not maximal and $f$ has an infinite number of saddle periodic points in 
$S(\C)\setminus S(\R).$ 
\end{thm}

\begin{proof} Let us assume that $f$ has maximal entropy and that $S(\R)$ 
has at least one bounded connected component $S(\R)^0.$
The existence of a bounded component implies that $S(\R)$ has five connected
components, four of which are unbounded and homeomorphic to disks, and 
one, $S(\R)^0,$ is bounded and homeomorphic to a sphere 
(see \cite{Benedetto-Goldman:1999}).
Being $f$-invariant and compact, $S(\R)^0$ is contained 
in $K(f).$ Since $f_\R$ has maximal entropy, $K(f)$ is 
contained in $S(\R).$ Since $K(f)$ is the support of $\mu_f$ 
(see \S \ref{par:henonlike}),
 $\mu_f(S(\R)^0)$ is a positive number. 
The ergodicity of $\mu_f$ and the $f$-invariance of $S(\R)^0$ now imply that $S(\R)^0$ 
has full $\mu_f$ measure. As a consequence, $K(f)$ coincides with $S(\R)^0.$ 
Since $S(\R)^0$ is compact, there is no one-sided periodic point, and
theorem \ref{thm:bscme} implies that $K(f)$ 
is a hyperbolic set. This means that the dynamics of $f$ on $S(\R)^0$ is 
uniformly hyperbolic. In particular, the unstable directions of $f$ 
determine a continuous line field on $S(\R)^0,$ and we get a  contradiction 
because $S(\R)^0$ is a sphere.  \end{proof}

\begin{cor}
Let $D$ be a real number and $S_D$  be the element of $\SS$ defined
by the real parameters $(0,0,0,D).$ The following properties are
equivalent: $(i)$ there exists a hyperbolic element $f$ in $\A$ such that
$f:S_D(\R)\to S_D(\R)$ has maximal entropy, $(ii)$
any hyperbolic element $f$ in $\A$ has maximal entropy
on $S_D(\R),$ and $(iii)$ $D\geq 4.$
\end{cor}

\begin{proof}
If $D>4,$ then $S(\R)$ is connected and smooth and the result follows
from theorem \ref{thm:entropy_minoration}. If $D\leq 0,$ the result 
follows from the fact that the action of the mapping class group on $S(\R)$ 
is totally discontinuous (see \cite{Goldman:2003}). If $0<D<4,$ then $S(\R)$ has a  compact connected component 
$S(\R)^0$ and the conclusion follows from the previous theorem.  \end{proof}


\subsection{Uniform hyperbolicity}\label{par:rduh}


We now prove theorem \ref{thm:hypreal} in the following more general form.

\begin{thm}\label{thm:hypreal+}
Let $f$ be a hyperbolic element of $\A.$ Let $S$ be an element of $\SS$
defined by real parameters. If $S(\R)$ is connected, then
\begin{itemize}
 \item the entropy of $f_\R$ is maximal; its value is $\log(\lambda(f));$
\item the set of bounded orbits of $f:S(\C)\to S(\C)$ is a compact subset $K(f)$ of  $S(\R);$
\item the automorphism $f$ admits a unique invariant probability measure $\mu_f$ of maximal entropy, and the support of $\mu_f$ coincides with $K(f)$; periodic saddle points equidistribute
toward $\mu_f;$
\item the dynamics of $f$ on $K(f)$ is uniformly hyperbolic.
\end{itemize} \end{thm}

The only property that has not been proven  yet is the last one. In fact, we shall prove more
than uniform hyperbolicity: Our objective includes a complete description of the complement of $K^+(f),$ in order to explain pictures like the one provided in figure \ref{fig:stables-num}. This will be 
achieved in section \ref{par:MacKay}.

\begin{figure}[t]
\input{ch-stable.pstex_t}
\caption{{\sf{Examples of stable manifolds.}} }
\label{fig:stables-num}
\end{figure}


\subsubsection{Notations and preliminaries}\label{par:prelimproof}
In what follows, we fix a hyperbolic  element $f$ of $\A,$  and assume
that $f$ preserves orientation (replace $f$  by $f^2$ if $f$ reverses orientation).
We denote by  ${\mathcal{H}}$   the space of real parameters $(A,B,C,D)$ such that $S(\R)$
is connected.
In order to prove theorem \ref{thm:hypreal+}, and theorem \ref{thm:mackay}, we shall
study the dynamics of $f$ on all surfaces $S=S_{(A,B,C,D)}$ with $(A,B,C,D)$ in ${\mathcal{H}}.$
For such surfaces,  maximal entropy implies the following properties:
\begin{enumerate}
\item $K(f)$ coincides with $J(f)$ and is a subset of $S(\R)$; moreover,
periodic points are hyperbolic, all of them are contained in $K(f)$, and intersections between stable and unstable manifolds are also contained in $K(f)$ 
(see theorem \ref{thm:bsrh}) ;
\item the set of one-sided points is a finite subset $OS(f)$ of $J(f)$ (see \cite{BS:2004}, 
sections 3 and 4)~; 
\item if $m$ is a point of tangency between stable and unstable manifolds
of $f,$ the $\alpha$ and $\omega$-limit sets of $m$ are contained in $OS(f)$
(see theorem 2.7 of \cite{BS:2004}) ~;
\item in the complement of $OS(f),$ stable and unstable manifolds of $f$ form
two laminations of $J(f)$ (see proposition 5.3 of \cite{BS:2002})~; 
\item a tangency  between a stable and an unstable manifold is always quadratic (see
section 2 in\cite{BS:2004}, section 5 of \cite{BS:2002}).
\end{enumerate}

Once again, as in the proof of theorem \ref{thm:entropy_minoration}, the main argument 
is to understand perturbations of the Cayley cubic, {\sl{i.e.}} perturbations of  
$f_\R:S_C(\R)\to S_C(\R).$


\subsubsection{Conical singularities of the Cayley cubic.}

The surface $S_C$ has
four conical singularities. If we blow up $\C^3$ at the four singular points of $S_C,$ 
the strict transform ${\widehat{S_C}}$ of $S_C$ is smooth, and the singular points are replaced
by four projective lines $\P^1(\C).$ These projective lines are called  exceptional 
divisors and will be denoted by $E_i,$ $i=1,$ $2,$ $3,$ $4.$
Another way to get the same surface $\widehat{S_C}$ is 
to blow up $\C^*\times \C^*$ at the four fixed points of the involution 
$
\eta(u,v)=(1/u,1/v).
$
The involution $\eta$ and the action of the group $\Gamma_2^*$ can be lifted to 
$\widehat{\C^*\times \C^*},$ and the quotient $\widehat{\C^*\times \C^*}/\eta,$ is 
isomorphic to $\widehat{S_C}.$ 

On the exceptional divisors, each hyperbolic element $f$ of $\A$ now has two hyperbolic 
fixed points, instead of one singular fixed point (this is the reason why we assume that
$f$ preserves orientation; if $f$ reverses orientation, then $f$ has a pair of periodic 
points of period $2$ on each exceptional divisor). 

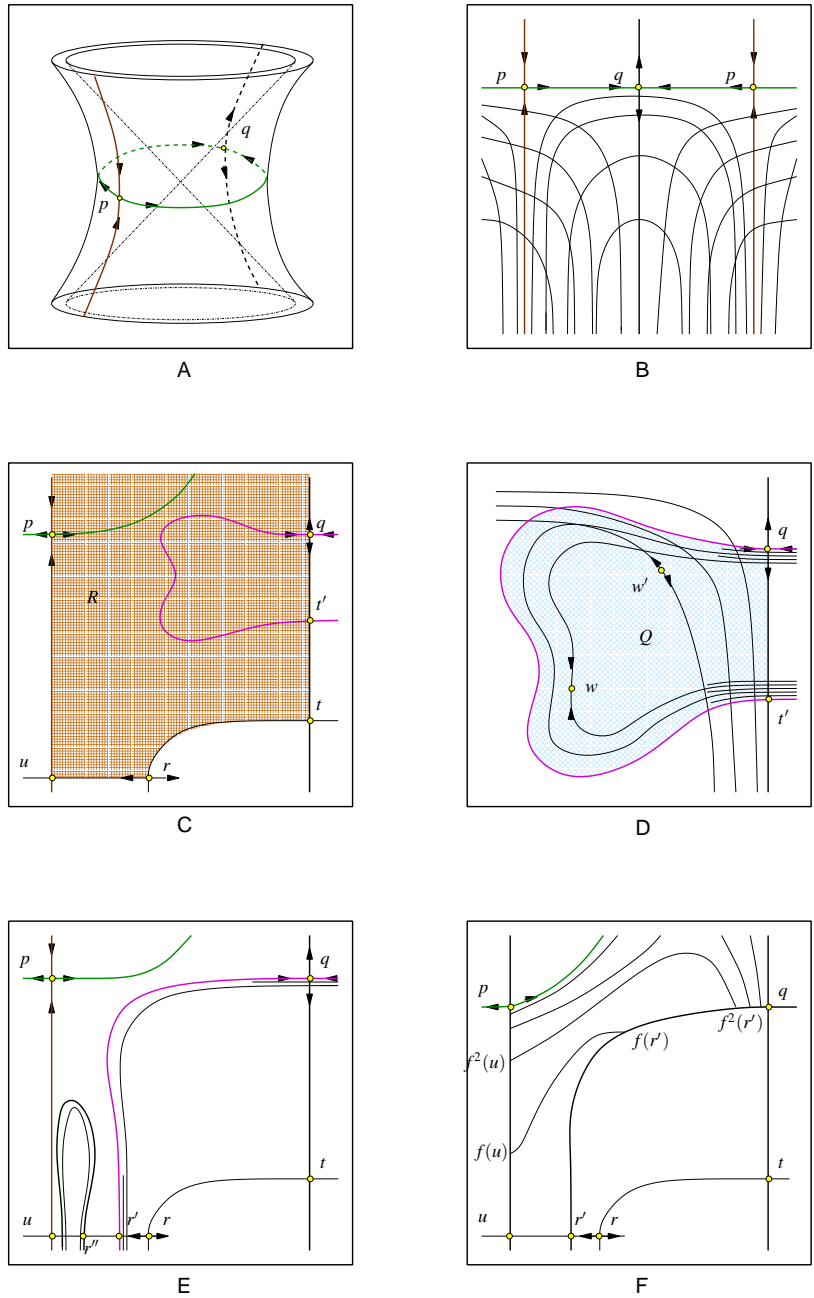
\begin{figure}[t]
\input{stables.pstex_t}
\caption{{\sf{Deformation of singularities.}} }
\label{fig:stables}
\end{figure}

Figures \ref{fig:stables} - {\sf{A}}  provides a local picture of $S_C(\R)$ after such a blow up.
Locally, this surface is a cylinder.   Figure \ref{fig:stables} - {\sf{B}}
is the same as figure \ref{fig:stables} - {\sf{A}}, but in the universal cover 
of the cylinder. It shows the geometry of the stable and unstable manifolds of
$f$ near $p$ and $q.$
These hyperbolic points are one-sided.  Moreover, the multiplier $\lambda(m)$ 
of $f$ along the unstable manifolds of these points are equal to $\lambda(f)^2,$ 
whereas multipliers of non singular periodic points
are equal to $\lambda(f).$ This illustrates a known property of one sided points (see
proposition 4.10 in \cite{BS:2002}). Last, but not least, 
the exceptional divisors are heteroclinic connections: Each exceptional 
divisor is at the same time the stable manifold of one of its periodic points and the 
unstable manifold of the other one. 

\begin{rem}\label{rem:nosc} The existence of such a heteroclinic connection is specific to this
construction: If $S$ is an element of $\SS$  and $f$ is a hyperbolic element
of $\A,$ there is no heteroclinic connection between periodic points of $f$ in $S(\C),$
because if $W^u(p)\setminus \{p\}$ coincides with $W^s(q)\setminus\{q\},$ then 
$W^u(p)\cup W^s(q)$  would be a one-dimensional compact complex subvariety of $\C^3.$
\end{rem}


\subsubsection{Deformation and periodic points.}\label{par:defpp}

For any point $(A,B,C,D)$ in ${\mathcal{H}},$ all periodic points of $f:S(\C)\to S(\C)$
are real saddle points (property (1) above). As a consequence, we can follow all the periodic
points along any deformation of the parameters $(A,B,C,D)$ in ${\mathcal{H}}$: 
If $\alpha(t),$ $t\in[0,1],$  is an arc of class ${\mathcal{C}}^k$ in ${\mathcal{H}},$ and if $p_0$ is a periodic saddle point of $f:S_{\alpha(0)}\to S_{\alpha(0)}$ 
of period $N,$  there exists an arc $p(t)$ of class ${\mathcal{C}}^k$  such that 
\begin{enumerate}
 \item for all $t,$ $p(t)$ is contained in $S_{\alpha(t)}$ and $p(0)=p_0$;
\item for all $t,$ $p(t)$ is a periodic saddle point of $f:S_{\alpha(t)}\to S_{\alpha(t)}$ of period $N$
(here we also use the fact that $f$ preserves orientation; otherwise, the period could change
when $p(t)$ goes through a singular point of $S_{\alpha(t)}$).
\end{enumerate}

\begin{rem}\label{rem:bjs}
The point $p(t)$ is contained in the set $K_{\alpha(t)}(f)$ of points in $S_{\alpha(t)}$ with
a bounded $f$-orbit. The family of compact sets  $K_{\alpha(t)}(f)$ depends semi-continu\-ously 
on $t$ (\cite{BS:1991inv}, lemma 3.1), so that the union 
$\cup_{t\in [0,1]} K_{\alpha(t)}(f)$
is contained in a fixed compact set ${\mathcal{K}}.$ The paths $p(t),$ 
$t\in [0,1],$ where $p$ describes the set of periodic points of 
$f:S_{\alpha(0)}\to S_{\alpha(0)}$ are all contained in 
${\mathcal{K}}.$
\end{rem}

Let us assume that $S_{\alpha(t)}$ is a deformation of the Cayley cubic $S_{\alpha(0)}=S_C.$ 
We shall say that a periodic point $w$ of 
$S_{\alpha(1)}$ {\sl{comes from a singular point $p$ of $S_C$}} if the deformation 
$w(t)$ of $w(1)=w$ along $\alpha(t)$ lands at  $p$ when $t=0.$

\begin{lem}\label{lem:one-sided-singular}
Let $m\in S$ be a periodic point of $f.$ The point $m$ is one-sided  if and only if
it comes from a singular point of the Cayley cubic. 
\end{lem}

\begin{proof}
Let us denote by $\alpha(t)$ a continuous path in the parameter space such that 
$S_{\alpha(0)}=S_C,$ $S_{\alpha(1)}=S$ and $\alpha(t)\in {\mathcal{H}},$ 
for all $t$ in $[0,1]$ ; such a path exists by remark \ref{rem:bgco} and theorem 
\ref{thm:entropy_minoration}.

Let $q$ be one of the four $u$-one-sided singular points of $S_C.$ Let us follow $q$ along the 
deformation and assume that $q(1)$ is not $u$-one-sided. We can then find a periodic 
point $w$ in $S$ such that $W^s(w)$ intersects $W^u(q(1))$ in two points $a$ and $b,$ 
one in each component of $W^u(q(1))\setminus\{ q(1)\}.$ 
We can now  follow $w$ and the intersection points $a$ and $b$  along the deformation
$S_{\alpha(t)}.$ If a tangency between $W^s(w(t))$ and $W^u(q(t))$ occurs at $a(t),$ 
the order of contact is $2$ (see section \ref{par:prelimproof}, property (5)).
After the tangency, there are two points of intersection 
between $W^s(w(t))$ and $W^u(q(t)),$ and both of them are contained in 
$S_{\alpha(t)}(\R)$  (see section \ref{par:prelimproof}, property (1)). 
Between these two points, we choose the
one which is closest to $q(t)$ along the component of 
$W^u(q(t))\setminus \{q(t)\}$ that contains it, and then continue to follow that point
along the deformation.
We apply the same strategy to $b(t).$

By semi-continuity   of the Julia set, $a(t)$ and $b(t)$ stay bounded when $t$ describes 
$[0,1]$ (see remark \ref{rem:bjs}). Since $q$ is $u$-one-sided, $a(0)$ and $b(0)$ 
are both on the same component of $W^u(q)\setminus \{q\}.$ From this we 
deduce the existence of a parameter $t_0$ such that $a(t_0)=q(t_0)$ 
(or $b(t_0)=q(t_0)$), which implies that there is a saddle connection between 
$w(t_0)$ and $q(t_0).$ This contradicts remark \ref{rem:nosc}.

This argument shows that $u$-one-sided (resp. $s$-one-sided) singular points of 
$S_C$ remain $u$-one-sided  (resp. $s$-one-sided) along any smooth deformation  
$\alpha(t)$ in ${\mathcal{H}}.$ The same argument shows that
the non-singular periodic points of $S_C$ cannot become one-sided. \end{proof}


\subsubsection{Deformation and stable manifolds}\label{par:exit}

Next steps aim at giving a description of $K(f)$ and are not absolutely necessary
to prove the uniform hyperbolicity. We shall make use of figures  \ref{fig:stables}-{\sf{B}}, to {\sf{E}};
they represent the geometry of stable and unstable manifolds near $p$ and $q$ after
deformation of the Cayley cubic.

Let us study the topology of  stable and unstable manifolds of $f$ on a  
connected deformation $S(\R)$ of $S_C(\R).$ For this, we consider one  exceptional 
divisor $E$ of ${\widehat{S_C}}$  and the two fixed points  $p$ and $q$ 
of $f$ on $E.$ Permuting $p$ and $q$ if necessary, we know that  $q$ is $u$-one
sided, half of its real unstable manifold going to infinity, and $p$ is $s$-one sided
(see picture \ref{fig:stables}-{\sf{B}}). 
We fix a periodic point $r$ in $S_C$ which is close to the stable manifold of $p$: 
The local unstable 
manifold of $r$ intersects transversaly the stable manifold of $p$ at $u$ and its stable 
manifold intersects transversaly the unstable manifold of $q$ at $t,$ as in 
 figure \ref{fig:stables}-{\sf{C}}. Changing $f$ in one of its iterates, we assume 
that  $r$ is a fixed point. 

We follow this picture along a small deformation $S_{\alpha(t)}$ between $S_C$ and 
$S=S_{\alpha(1)},$ keeping the same local geometry for $W^s(p),$ $W^u(r),$ $W^s(r),$ and $W^u(q).$

Let $R\subset S(\R)$ be the closed region which is bounded by  
the half of  $W^s(p)\setminus \{u\}$ that contains $u,$
the segment of $W^u(r)$ between $u$ and $r,$
the segment  of $W^s(r)$ that joins $r$ to $t,$ and
the half of $W^u(q)\setminus \{t\}$ that contains $t,$ 
(see figure \ref{fig:stables}-{\sf{C}}).
Let $W^s_+(q)$ be  the connected component of $W^s(q)\setminus \{q\}$ which enters $R$:
this half stable manifold is parametrized by $\xi:\R^+\to S(\R),$ with $\xi(0)=q$ and $\xi(z)\in R$
for small positive real numbers $z.$  
The closure of the stable manifold of $q$ covers the set $K(f).$ As a consequence, 
we may assume that $W^s_+(q)\setminus \{q\}$ exits the region $R$ (in particular, 
there exists a positive $z$ such that $\xi(z)$ is on $\partial R$).

\begin{lem}\label{lem:exit}
The half stable manifold $W^s_+(q)$ exits $R$ through $W^u(r),$ in between $r$ and $u.$ 
\end{lem}

\begin{proof}
The stable manifold $W^s(q)$ must exit  $R$ through a piece of unstable manifold. 
We may therefore suppose that it leaves $R$ through $W^u(q),$ in between $q$ and the point $t$ (see picture \ref{fig:stables} - {\sf{C}} and {\sf{D}}). Let $t'$ be the first point 
of intersection 
of $W^s(q)$ and the segment $[q,t]\subset W^u(q),$ and let $Q$ be the subset 
of $R$ which is bounded by the arcs of  $W^s_+(q)$ and $W^u(q)$ in between 
$q$ and $t.$ 

Since $q$ is $u$-one-sided, and since the stable manifolds form a lamination near $t'$ 
(see section \ref{par:prelimproof}, property (4)), we know that there is no stable manifold that
approaches $t'$ along the segment $[t',t]\subset W^u(q).$ 
On the other hand, $t'$ is not isolated in $K(f),$ and, in particular, there are stable manifolds
 that cross $W^u(q)$ in between $q$ and $t',$ and arbitrarily near $t'.$ Those manifolds 
 "follow" $W^s_+(q)$: They enter $Q$ near $t'$ and then exit $Q$ near $q$ (see picture 
 \ref{fig:stables} - {\sf{D}}). 
 Let now $W$ be  an unstable manifold which intersects $W^s_+(q)$ in between $q$ 
 and $t'$: Such a manifold enters $Q,$ and then has to leave it, intersecting 
 $W^s_+(q)$ a second time between $q$ and $t'.$ These two observations and 
 the density of periodic points in $K(f)$ imply that $Q$ contains periodic points 
 $w$ and $w'$ of $f$ such that the connected
 component of $W^u(w)\cap Q$ containing $w$ intersects the connected 
 component of $W^s(w')\cap Q$ containing $w'$ in at least two distinct points.

 We now follow this picture along a deformation $\alpha(\epsilon)$ between 
 $S=S_{\alpha(1)}$
 and the Cayley cubic $S_C=S_\alpha(0),$ as in lemma \ref{lem:one-sided-singular}. 
 Periodic points, stable and unstable manifolds,  and intersections between these curves 
 move continuously along the deformation.
  In particular, 
 the periodic points $w(\epsilon)$ and $w'(\epsilon),$ and the (at least) two intersection points
 of their stable/unstable manifolds stay in the region $R(\epsilon).$ 
 Let $N$ be a common 
period for $w$ and $w'.$ 
The number of points in the region $R(0)$ on the Cayley cubic 
with period $N$ is finite, all of them are periodic saddle points 
(except for $p$ and $q$) and the local stable and unstable manifolds of these points 
in the region $R(0)$ intersect in exactly one point (see the proof of lemma 
\ref{lem:one-sided-singular} and picture \ref{fig:stables} - {\sf{B}}). 
From this we deduce the existence of a parameter $\epsilon\geq
 0$ for which $W^u(w(\epsilon))\cap W^s(w'(\epsilon))$ is not contain in $S_{\alpha(\epsilon)}(\R)$:
At least one intersection points has become complex.
Once again, this contradicts the fact that $f_\R:S_{\alpha(\epsilon)}(\R)\to S_{\alpha(\epsilon)}(\R)$ 
has maximal entropy (see section \ref{par:prelimproof}, property (1)).

This shows that $W^s(q)$ cannot leave the region $R$ through $W^u(q).$ The only remaining 
possibility is  that $W^s(q)$ leaves $R$ through $W^u(r),$ in between $r$ and $u.$ 
\end{proof}


\subsubsection{Deformation, stable manifolds and doubly one-sided points}\label{par:defsmdo}

Let $I$ be the closed segment $[r,u]\subset W^u(r).$
Let $r'$ be the first point of intersection of $W^s_+(q)$ with 
$I$ (see figure \ref{fig:stables}-{\sf{E}}). 
Since $q$ is $u$-one-sided, we know that there is no stable manifolds  
approaching $r'$ from the left. 
We can therefore define $r''$ to be the unique point in $I$ which is between $u$ and $r',$  
]is contained in $K(f),$ and  is closest to $r'$ with these properties.

If $r''$ is different from $u,$ the stable manifold through $r''$ enters $R$ and cannot intersect
$W^s(q)$ and $W^s(p).$ It must therefore exit $R$ through the interval $I,$ in between 
$r''$ and $u$ (see picture \ref{fig:stables} - {\sf{E}}). We then obtain a contradiction along 
the same line as in lemma \ref{lem:exit}. This implies that 
$r''$ coincides with the point $u,$ and no stable manifold crosses $I$ between $r'$ and $u.$  

\begin{figure}[t]\label{fig:strip}
\input{strip.pstex_t}
\caption{ {\sf{Complement of $K^+(f).$}} }
\end{figure}
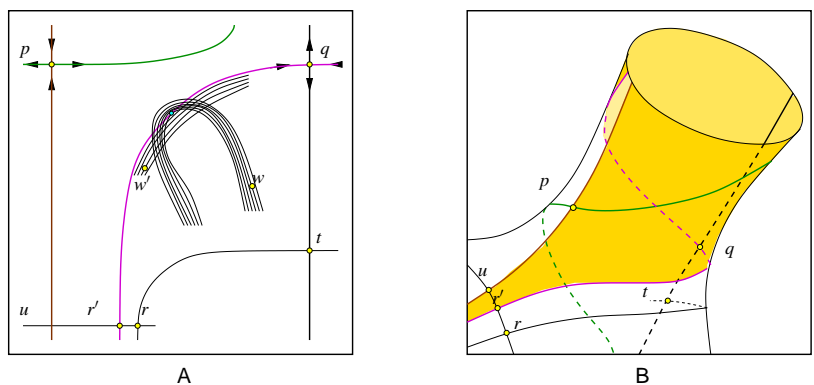

The segments $f^n[r',u]$ join the endpoints $f^n(u),$ which converge 
to $p$ along $W^s(p),$ to $f^n(r'),$ which converge to $q$ along $W^s(q).$ 
These segments are pieces of unstable manifolds and, as such, they can not 
intersect $W^u(p).$ This shows that the connected component of $W^u(p)$ 
that enters $R$ does not leave $R$: This half unstable manifold has to go 
to infinity, and $p$ is both $u$ and $s$-one-sided
(see picture \ref{fig:stables} - {\sf{F}}) .


\subsubsection{Deformation and the geometry of $K(f)$}\label{par:defgeo}

We can apply the same argument to understand the geometry of stable and unstable 
manifolds near  $p.$ Part {\sf{B}} of figure \ref{fig:strip} summarizes our knowledge of the 
geometry of stable and unstable manifolds near the points $p$ and $q$ after a small 
deformation of the Cayley cubic: $p$ and  $q$ are both $u$ and $s$-one-sided, 
and the colored region is contained  in the complement of $K(f).$ 

Let us now consider a large deformation $S_{\alpha(t)}$ of the 
Cayley cubic $S_C.$ Following $p,$ $u,$ $r,$ $t,$ $q,$ and the
stable/unstable manifolds of these points along the deformation, 
we can follow the region $R$ along $\alpha(t).$ Since there
is no saddle connection in $S_{\alpha(t)}$ for $t\neq 0,$ the geometry 
of $R$ with respect to local stable and unstable manifolds in $R$ does not change.
The results obtained above for small deformation remain therefore valid for 
arbitrarily large deformation $\alpha(t) \subset {\mathcal{H}}.$


\subsubsection{Absence of tangency and hyperbolicity}
Let us assume 
that there is at least one set of parameters $(A,B,C,D),$ for which $S(\R)$ is connected
and $f_\R$ is not uniformly hyperbolic along $K(f).$ Then, there is a tangency 
between the stable manifold
of a $u$-one-sided periodic point $q$ and an unstable manifold. Iterating $f,$ we can find such tangencies
in arbitrarily small neighborhoods of $q.$ Since $q$ is $u$-one-sided, the previous 
steps describe the geometry of the stable and unstable manifolds near $q.$ 
Figure \ref{fig:strip}-{\sf{A}} represents such a possible 
tangency, and it shows the geometry of the unstable lamination near such a tangency 
(see \cite{BS:2004}, picture 4.1  and sections 3 and 4, for a detailed description). 

In an arbitrarily small
neighborhood ${\mathcal{U}}$ of the tangency, we can find a periodic 
saddle point $w,$ such that the connected 
component $W^u_{loc}(w)$ of $W^u(w)\cap {\mathcal{U}}$ containing $w$ 
intersects $W^s(q)$  in two points. Then, we can find a second periodic 
saddle point $w'$ such that $W^s(w')$ intersects $W^u_{loc}(w)$
in two points (see figure \ref{fig:strip}-{\sf{A}}). 

Since $S(\R)$ is connected, we can follow this picture along a deformation 
$\alpha(\epsilon)\in {\mathcal{H}}$, $\epsilon\in (0,1],$ which approaches the Cayley parameters $(0,0,0,4)$
when $\epsilon$ goes to $0$ (see remark \ref{rem:bgco}). We know from section \ref{par:defpp}
that the periodic saddle points $r(\epsilon),$ $t(\epsilon),$ $u(\epsilon),$ $w(\epsilon)$ and $w'(\epsilon),$
move continuously along the deformation. 
From sections \ref{par:defsmdo} and \ref{par:defgeo} , we can also assume that  
the geometry of the stable and unstable manifolds of $r(\epsilon),$  $t(\epsilon),$ $q(\epsilon)$ and 
$p(\epsilon)$ remains the same  along the deformation; in particular,
since the periodic points $w(\epsilon)$ and $w'(\epsilon)$ cannot cross the stable or 
unstable manifolds of other periodic points during the deformation,  they both 
stay in the interior of the region $R(\epsilon).$  We then get a contradiction as in section \ref{par:defsmdo}.

Since there is no tangency, the dynamics of $f$ is uniformly hyperbolic on $K(f).$ This proves theorem \ref{thm:hypreal+}.


\subsection{Strips and bounded orbits}\label{par:MacKay}


Let $(A,B,C,D)$ be an element of ${\mathcal{H}}.$
Let $f$ be a hyperbolic  element of $\A.$ 
The surface $S(\R)$ defined by this
set of parameters is connected, and $f:S(\R)\to S(\R)$ is uniformly hyperbolic on $K(f),$
so that we can apply proposition 2.1.1 of \cite{Bonatti-Langevin:ast}: The set
$$
W^s_\R(K(f))=K^+(f)\cap S(\R)
$$
is laminated by stable manifolds of points in $K(f)$; if a point $m$ in $K^+(f)$ is on 
the boundary of the complement of $W^s_\R(K(f)),$ then $m$ is on the stable 
manifold of
a periodic $u$-one-sided periodic point of $f.$ From section \ref{par:defpp}, we know
that $f$ has exactly eight periodic one-sided points, each of them coming from a 
singularity of the Cayley cubic. From sections \ref{par:exit} and \ref{par:defgeo}, the
stable manifolds of the two one-sided points coming from one singularity bound a strip,
as in picture \ref{fig:strip}-{\sf{B}}.
This proves the following result, which was first numerically observed by Catarino and MacKay (see \cite{Catarino:2004}, page 61 for example), and "explains" pictures \ref{fig:stables-num}-{\sf{A,C}}.

\begin{thm}[MacKay observation]\label{thm:mackay}
If $S(\R)$ is connected, $f$ has exactly eight one-sided  fixed points $p_1,$ $q_1,$
$p_2,$ $q_2,$   $p_3,$ $q_3,$ $p_4,$ and $q_4.$ All of them come from singularities of
the Cayley cubic by deformation; all of them are both $u$ and $s$-one-sided.
Moreover, the stable manifolds of $p_i$ and $q_i$ ($i=1,2,3,4$) bound an open strip homeomorphic 
to $\R\times (-1,1),$ and the complement
of $K^+(f)\cap S(\R)$ coincides with the union of these four strips. \end{thm}

\begin{rem}
We shall prove in theorem \ref{thm:hausdorff} that the Hausdorff dimension of $K^+(f)\cap S(\R)$
is stricly less than $2.$ In particular, its complement, {\sl{i.e.}} the union of the four strips, has
full Lebesgue measure ; almost all orbits go to infinity under iteration of $f.$ The same
is true for the complement of $K^+(f)$ in $S(\C).$  
\end{rem}


\section{Schr\"odinger operators and Painlev\'e equations}\label{par:schrodinger}

\subsection{Discrete Schr\"odinger operators} Let us now apply the previous results to 
the study of the spectrum of certain discrete 
Schr\"odinger operators. There is a huge literature on the subject, and we refer to  
\cite{Damanik:2000} and \cite{Damanik:survey} for background results and a short
bibliography.


\subsubsection{Discrete Schr\"odinger operators and substitutions}
Let $W^*$ be the set of finite words in the letters $a$ and $b.$ Let $\iota:\{a,b\}
\to W^*\setminus \{ \emptyset \}$ be a substitution.  In what follows, we will assume that
$\iota$ is invertible, which means that $\iota$ extends to an automorphism
$\Phi_\iota$ of the free group $F_2=\langle a, \, b \, \vert \, \emptyset \, \rangle,$
and that $\iota$ is primitive, which means that $\Phi_\iota$ is hyperbolic ; 
in other words, the image of $\Phi_\iota$ in $\Out(F_2)=\GL(2,\Z)$ 
is a hyperbolic matrix, with two distinct eigenvalues $\lambda_+(\iota)$
and $\lambda_-(\iota)$  satisfying
$$
\vert \lambda_+(\iota) \vert = \vert 1/\lambda_-(\iota)\vert >1. 
$$
Under these hypothesises, there is a unique infinite word $u_+$ in the two letters
$a$ and $b$ such that  $\iota(u_+)=u_+.$

\begin{eg} The Fibonacci substitution $\iota_F,$ defined by
$\iota_F(a)=b$ and $\iota_F(b)=ba,$ provides a good and famous example of such 
an invertible primitive substitution. Its fixed word starts with
$babbababbabbababbababbabba...$
\end{eg}

Let $W$ be the set of bi-infinite words in $a$ and $b$ and $\tilde{T}:W\to W$
be the left shift. Let $\tilde{u}_+$ be any completion of $u_+$ on the left. We then 
define $\Omega$ to be the $\omega$-limit set of the $\tilde{T}$-orbit of $\tilde{u}_+$:
$$
\Omega = \left\{
v \in W\, \vert \quad  {\text{there exists a sequence }} n_i\to +\infty, {\text{ such that }} {\tilde{T}}^{n_i}({\tilde{u}}_+)\to v
\right\} .
$$
Since $\iota$ is primitive, the restriction of the left shift 
${\tilde{T}}$ to the set $\Omega$ is a minimal and uniquely ergodic
homeomorphism $T:\Omega\to \Omega.$ The unique $T$-invariant 
probability measure on $\Omega$ will be denoted by $\nu.$

\begin{rem} The subshift $T:\Omega \to \Omega $ encodes the dynamics
of a rotation $R_\alpha : \R/\Z \to \R/\Z,$ where $\alpha$ is 
a quadratic integer (see \cite{Arnoux:2002}).This 
provides a measurable conjugation between $R_\alpha$ and $T$ 
which sends the Lebesgue measure $dx$ to $\nu.$
\end{rem}

Let us now fix an element $w$ in $\Omega,$ and define the potential 
$V_w:\Z\to \R$ by  $V_w(n)=1$ if $w_n=a$ and $V_w(n)=0$ if $w_n=b.$
Let $\kappa$ be any complex number ($\kappa$ is the so called 
"coupling parameter"). If $(\xi(n))_{n\in \Z}$ is a complex valued sequence, 
we define 
$$
H_{\kappa, w}(\xi)(n)= \xi(n+1)+\xi(n-1)+\kappa V_w(n)\xi(n).
$$ 
The discrete Schr\"odinger operator $H_{\kappa, w}$ induces a bounded linear
operator on $l^2(\Z),$ with norm at most $2+\vert \kappa \vert.$ The 
adjoint of $H_{\kappa, w}$ is $H_{{\overline{\kappa}}, w},$ so that $H_{\kappa, w}$
is self-adjoint if and only if $\kappa$ is a real number. 


\subsubsection{Almost sure spectrum and Lyapunov exponent}
Since $T$ is ergodic with respect to $\nu,$ there exists a subset $\Sigma_\kappa$ 
of $\C$ (of $\R$ if $\kappa$ is real) such that the spectrum of 
$H_{\kappa,w}:l^2(\Z)\to l^2(\Z)$ coincides with $\Sigma_\kappa$ for $\nu$-almost all $w$ 
in $\Omega.$ This set is the "almost sure spectrum" of the family $H_{\kappa, w}.$ 

To understand the spectrum of $H_{\kappa, w},$ one is led to solve the eigenvalue 
equation $H_{\kappa, w}(\xi) = E \xi$ ($E$ in $\R$ or $\C$). For any initial condition 
$(\xi(0),\xi(1)),$ there is a unique solution, which is given by the recursion formula
$$
\left(\! \begin{array}{c} \xi(n+1) \\ \xi(n) \end{array}\!  \right) = \left(\!  \begin{array}{cc} E-\kappa V_w(n) & -1 \\ 1 & 0 \end{array} \! \right) 
\left(\!  \begin{array}{c} \xi(n) \\ \xi(n-1) \end{array}\!  \right), \quad n\in \Z.
$$
Let $M_{\kappa,E}:W^*\to \SL(2,\C)$ be defined by 
$$
M_{\kappa, E} (a)= \left( \begin{array}{cc} E-\kappa & -1 \\ 1 & 0 \end{array}\right), \quad 
M_{\kappa, E} (b)= \left( \begin{array}{cc} E & -1 \\ 1 & 0 \end{array}\right), 
$$
and by 
$$
M_{\kappa,E}(u_1...u_n)=\Pi_{i=0}^{i=n-1} M_{\kappa, E}(u_{n-i})
$$
for any word $u=u_1...u_n$ of length $n.$ 
This defines  a $\SL(2,\C)$-valued cocyle over the dynamical system $(\Omega, T, \nu).$ Applying 
Osseledets' theorem, each choice of a coupling parameter $\kappa$ and an energy $E$ gives rise 
to a non negative Lyapunov exponent $\gamma(\kappa, E),$ such that
\begin{eqnarray*}
\gamma(\kappa,E) & = & \lim_{n\to + \infty} \frac{1}{n} \int_\Omega \log \Vert M_{\kappa, E}(w_1 w_2...w_{n-1})\Vert \, d\nu(w) \\
  & = & \lim_{n\to +\infty} \frac{1}{n} \log  \Vert M_{\kappa, E}(w_1 w_2...w_{n-1})\Vert,
\end{eqnarray*}
for $\nu$-almost all $w$ in $\Omega.$ The Lyapunov function $\gamma(\kappa,E)$ is linked to the
almost sure spectrum $\Sigma_\kappa$ by the following result. 

\begin{thm}[see \cite{Damanik:2000}]
Let $\kappa$ be a real number. The almost sure spectrum $\Sigma_{\kappa}$ coincides with the set 
of energies for which the Lyapunov exponent vanishes.
\end{thm}


\subsubsection{Trace map dynamics, Lyapunov exponent, and Hausdorff dimension}
Let us fix the coupling parameter $\kappa.$ Let $S_{4+\kappa^2}$ be the character surface 
$x^2+y^2+z^2-xyz=4+\kappa^2.$
The {\sl{Schr\"odinger curve}} of $S_{4+\kappa^2}$ is the parametrized rational curve
$s:\C\to S_{4+\kappa^2},$ which is  defined by $s(E)=(x(E),y(E),z(E)),$ with 
\begin{eqnarray*}
(x(E),y(E),z(E)) & = & (\tr(M_{\kappa,E}(a)),\, \tr(M_{\kappa,E}(b)),\, \tr(M_{\kappa,E}(ab))) \\
                         & = & (E-\kappa,\, E ,\, E(E-\kappa)-2).
\end{eqnarray*}

\begin{rem}
The intersection of $S_{4+\kappa^2}$ with the plane $y=x+\kappa$ is a reducible cubic curve: It is the
union of $s(\C)$ with the line $\{z=2, \, y=x+\kappa\}$; the involution 
$s_z$ permutes these two curves.
\end{rem}

Let $f_\iota$ be the polynomial automorphism of $S_{4+\kappa^2}$ which is determined by the automorphism
$(\Phi_\iota)^{-1}:F_2\to F_2.$ By definition of $f_\iota,$ we have 
$$
\left(\tr(M_{\kappa,E}(\iota(a))),\, \tr(M_{\kappa,E}(\iota(b))),\,  \tr(M_{\kappa,E}(\iota(ab)))\right)= f_\iota (s(E)).
$$
In \cite{Damanik:2000}, Damanik proved that $\gamma(\kappa,E)$ vanishes if and only if the point $s(E)$ 
has a bounded forward $f_\iota$-orbit. In other words, $\Sigma_\kappa$ is given by the intersection 
between the Schr\"odinger curve and the set $K^+(f_\iota)$: 
\begin{equation}\label{eq:spectrum}
\Sigma_\kappa=  \left\{  E\in \C\, \vert \, s(E)\in K^+(f_\iota)  \right\}.
\end{equation}

We can now apply MacKay observation, i.e. theorem \ref{thm:mackay}, which tells
us that the complement of $s(\Sigma_\kappa)$ in the real Schr\"odinger curve is  
obtained by intersecting $s(\R)$ with the four strips associated to the one-sided
points of $f.$ This means that {\sl{gaps in the complement of the spectrum are bounded
by intersection points between $s(\R)$ and the eight curves 
$W^s(q_i)$  and $W^s(p_i),$ $i=1,$ $2,$ $3,$ and $4$)}}.

\begin{thm}\label{thm:hausdorff}
The Hausdorff dimension of $\Sigma_\kappa,$ $\kappa\in \R,$ is a real analytic function of $\kappa^2.$ Moreover, 
$$
0< Haus(\Sigma_\kappa) \leq 1, \quad \forall \kappa \in \R,
$$ 
and $Haus(\Sigma_\kappa) = 1$ if and only if $\kappa = 0.$
\end{thm}

This statement confirms numerical observations that can be found, for example, 
in  \cite{Kohmoto-Kadanoff-Tang:1983} and \cite{Kohmoto:1983}; it is stronger than the fact that $\Sigma_\kappa$ has zero Lebesgue measure 
when $\kappa\neq 0,$ a property which was proved by Kotani in the eighties 
(see \cite{Damanik:survey}). Here, it appears as a corollary of results in dynamical
systems which are due to Bowen, Pesin, and Ruelle. 

\begin{proof}
When $\kappa$ is a non zero real number, we obviously have $4+\kappa^2>4,$ and theorem 
\ref{thm:hypreal} shows that the dynamics of $f_\iota$ is 
uniformly hyperbolic on its Julia set. 
By results of Hasselblatt \cite{Hasselblatt:1994}, 
the stable and unstable distributions of $f_\iota$ are smooth, and the holonomy maps between two transversals
of the stable (resp. unstable) laminations are Lipschitz continuous. In particular, the Hausdorff dimension 
of the sets 
$$
W^u_{loc}(m)\cap K^+(f_\iota)
$$
does not depend on the choice of $m$ in $K(f_\iota),$ and, by (\ref{eq:spectrum}), it coincides with the dimension of $\Sigma_\kappa.$
The map $f_\iota$ is area-preserving:
As in \cite{Wolf:2002}, corollary 4.7, this implies that the Hausdorff dimension of the 
sets $W^s_{loc}(m)\cap K^-(f_\iota)$ coincides also with $Haus(\Sigma_\kappa).$

Using Bowen-Ruelle thermodynamic formalism,
as it is done in \cite{Verjovsky-Wu:1996}, we obtain that the Hausdorff dimension of $\Sigma_\kappa$ is an analytic
function of $\kappa^2.$ Since the function $G^+_{f_\iota \vert s(E)}$ is H\"older continuous
 the Hausdorff dimension of $\Sigma_\kappa$ is strictly positive. 

Let us now show that $Haus( \Sigma_\kappa)$ is strictly less than $1.$ If 
$Haus(\Sigma_\kappa)=1,$ the Hausdorff dimension of the slices 
$W^u_{loc}(m)\cap K^+(f_\iota)$ and $W^s_{loc}(m)\cap K^-(f_\iota)$ 
are also equal to $1.$ Theorem 22.1 of \cite{Pesin:book} then shows that the 
Lebesgue measure of these sets is strictly positive. By Hasselblatt's
result, the Lebesgue measure of $K(f_\iota)$ is positive, and by Bowen-Ruelle's 
theorem (\cite{Bowen-Ruelle:1975}, theorem 5.6), the set $K(f_\iota)$ must be an attractor 
of $f_\iota:S_{4+\kappa^2}\to S_{4+\kappa^2}$ This contradicts the
fact that $K(f_\iota)$ is compact and $f$ is area preserving.
\end{proof}

\begin{rem} It would be interesting to settle a complete dictionary between dynamics of
the trace map and properties of the spectrum. For example,  the Green function
of $f_\iota$ should coincide with the Lyapunov function $\gamma(\kappa,E)$ along the Schr\"odinger
curve ; together with Thouless formula, this would identify the density of states $dk_\kappa$ with the measure obtained by slicing $T^+_{f_\iota}$ with the Schr\"odinger curve: 
$dk_\kappa = s^*(T^+_{f_\iota})$ (see \cite{Sabot:2001} for related results, the definition of $dk_\kappa,$
and Thouless formula).
\end{rem}


\subsection{Monodromy of Painlev\'e VI equation}\label{par:painleve6}

The sixth Painlev\'e
 equation $P_{VI}=P_{VI}(\theta_\alpha,\theta_\beta,\theta_\gamma,\theta_\delta)$
is the second order non linear ordinary differential equation
\begin{equation*}
P_{VI}\ \ \ \left\{
\begin{matrix}{ {d^2 q\over dt^2}}&=&{{1\over 2}\left({1\over q}+{1\over q-1}+{1\over q-t}\right)
\left({dq\over dt}\right)^2
-\left({1\over t}+{1\over t-1}+{1\over q-t}\right)
\left({dq\over dt}\right)}\\
&&{+{q(q-1)(q-t)\over t^2(t-1)^2}
\left(\frac{\theta_\delta^2}{2}-\frac{\theta_\alpha^2}{2}{t\over q^2}
+\frac{\theta_\beta^2}{2}{t-1\over (q-1)^2}+\frac{1-\theta_\gamma^2}{2}{t(t-1)\over
(q-t)^2}\right)}.
\end{matrix}\right.
\end{equation*}
the coefficients of which depend on complex parameters $\theta=(\theta_\alpha,\theta_\beta,\theta_\gamma,\theta_\delta)$.

As explained in \cite{Iwasaki-Uehara:2006} (see also
\cite{Cantat-Loray:2007}), the monodromy of Painlev\'e equation provides a representation
of $\pi_1(\P^1\setminus\{0,1,\infty\},t_0)$ into the group of
analytic diffeomorphisms of the space of initial conditions $(q(t_0),q'(t_0))$ (see 
\cite{Iwasaki-Uehara:2006} for a precise necessary description of this space).
Via the Riemann-Hilbert correspondence, the space of initial conditions
is analytically isomorphic  to (a desingularization of) $S_{(A,B,C,D)}$ with 
parameters
\begin{equation}\label{theta/abcd}
a=2\cos(\pi\theta_\alpha),\ b= 2\cos(\pi\theta_\beta),\  c= 2\cos(\pi\theta_\gamma),\  d= 2\cos(\pi\theta_\delta),
\end{equation}
the monodromy action on  the space of initial conditions is conjugate 
to the action of $\Gamma_2$ on  the surface $S_{(A,B,C,D)}.$  

From this, and from sections \ref{par:rduh} and \ref{par:schrodinger}, we deduce the following result, thereby answering a recent question raised by Iwasaki and Uehara, as
problem 15 of \cite{Iwasaki:NagoyaQuestion}.

\begin{thm}\label{thm:painlevesingulier}
Let $(\theta_\alpha,\theta_\beta,\theta_\gamma,\theta_\delta)$ be parameters of Painlev\'e
sixth equation such that  
\begin{itemize}
\item[(i)] for $\epsilon= \alpha,$ $\beta,$ $\gamma,$ and $\delta,$ the real part of
$\Theta_\epsilon$ is an integer $n_\epsilon,$ and
\item[(ii)] $n_\alpha + n_\beta + n_\gamma+n_\delta$ is odd.
\end{itemize}
Let $\eta$ be any loop in $\P^1\setminus\{0,1,\infty\},$ and let 
$
f_\eta:S_{(A,B,C,D)}\to S_{(A,B,C,D)}
$ 
be the mondromy transformation defined by $\eta$ (through Riemann-Hilbert correspondence). 
If the entropy of $f_\eta$ is positive, then 
\begin{itemize}
\item all periodic points of $f_\eta$ are contained in the real part $S_{(A,B,C,D)}(\R)$
of the surface;
\item the Hausdorff-dimension of the set of bounded  $f_\eta$-orbits is less than
$2$;
\item the unique invariant probability measure of maximal entropy $\mu_{f_\eta}$ is
singular with respect to the Lebesgue measure on $S_{(A,B,C,D)}(\R)$ (and $S_{(A,B,C,D)}(\C)$).
\end{itemize}
\end{thm}

This theorem should be  compared to Goldman's results regarding  ergodic properties
of the whole $\Gamma_2^*$ action with respect to the invariant area form $\Omega$
(see the definition of $\Omega$ in section \ref{par:prelim-amg}). As a particular case of
Goldman's theorem, the action of $\Gamma_2^*$ on $S_D(\R)$ is ergodic with respect
to $\Omega$ if, and only if $4< D \leq 20$ (see  \cite{Goldman:2003}).
Another interesting example is given by the Markoff surface $S_0.$ In this example, the
quasifuchsian space $\QF$ provides an open invariant subset of $S_0(\C).$ This shows that
the action of $\Gamma_2$ on $S_0(\C)$ is not ergodic. Theorem \ref{thm:painlevesingulier}
and these results suggest that, for most parameters, the dynamics of the monodromy of
Painlev\'e equation is not correctly described by the invariant area form $\Omega.$ 

%
%

\bibliographystyle{plain}
\bibliography{referencesbhps}

\end{document}

%% file: intro-julia.pstex_t
\begin{picture}(0,0)%
\includegraphics{intro-julia.pstex}%
\end{picture}%
\setlength{\unitlength}{1776sp}%
\begingroup\makeatletter\ifx\SetFigFont\undefined%
\gdef\SetFigFont#1#2#3#4#5{%
  \reset@font\fontsize{#1}{#2pt}%
  \fontfamily{#3}\fontseries{#4}\fontshape{#5}%
  \selectfont}%
\fi\endgroup%
\begin{picture}(12644,6044)(1179,-6383)
\end{picture}%

%% file: ch-stable.pstex_t
\begin{picture}(0,0)%
\includegraphics{ch-stable.pstex}%
\end{picture}%
\setlength{\unitlength}{1579sp}%
\begingroup\makeatletter\ifx\SetFigFont\undefined%
\gdef\SetFigFont#1#2#3#4#5{%
  \reset@font\fontsize{#1}{#2pt}%
  \fontfamily{#3}\fontseries{#4}\fontshape{#5}%
  \selectfont}%
\fi\endgroup%
\begin{picture}(13244,20287)(1179,-20626)
\put(3901,-6661){\makebox(0,0)[lb]{\smash{{\SetFigFont{9}{10.8}{\familydefault}{\mddefault}{\updefault}{\sf{A}}-1}}}}
\put(11101,-6736){\makebox(0,0)[lb]{\smash{{\SetFigFont{9}{10.8}{\familydefault}{\mddefault}{\updefault}{\sf{A}}-2}}}}
\put(11101,-13636){\makebox(0,0)[lb]{\smash{{\SetFigFont{9}{10.8}{\familydefault}{\mddefault}{\updefault}{\sf{B}}-2}}}}
\put(3826,-13636){\makebox(0,0)[lb]{\smash{{\SetFigFont{9}{10.8}{\familydefault}{\mddefault}{\updefault}{\sf{B}}-1}}}}
\put(3901,-20536){\makebox(0,0)[lb]{\smash{{\SetFigFont{9}{10.8}{\familydefault}{\mddefault}{\updefault}{\sf{C}}-1}}}}
\put(11101,-20536){\makebox(0,0)[lb]{\smash{{\SetFigFont{9}{10.8}{\familydefault}{\mddefault}{\updefault}{\sf{C}}-2}}}}
\end{picture}%

%% file: stables.pstex_t
\begin{picture}(0,0)%
\includegraphics{stables.pstex}%
\end{picture}%
\setlength{\unitlength}{1184sp}%
\begingroup\makeatletter\ifx\SetFigFont\undefined%
\gdef\SetFigFont#1#2#3#4#5{%
  \reset@font\fontsize{#1}{#2pt}%
  \fontfamily{#3}\fontseries{#4}\fontshape{#5}%
  \selectfont}%
\fi\endgroup%
\begin{picture}(16866,27141)(1168,-27469)
\put(3076,-4636){\makebox(0,0)[lb]{\smash{{\SetFigFont{7}{8.4}{\familydefault}{\mddefault}{\updefault}$p$}}}}
\put(6076,-3061){\makebox(0,0)[lb]{\smash{{\SetFigFont{7}{8.4}{\rmdefault}{\mddefault}{\updefault}$q$}}}}
\put(13876,-1936){\makebox(0,0)[lb]{\smash{{\SetFigFont{7}{8.4}{\rmdefault}{\mddefault}{\updefault}$q$}}}}
\put(11401,-1936){\makebox(0,0)[lb]{\smash{{\SetFigFont{7}{8.4}{\familydefault}{\mddefault}{\updefault}$p$}}}}
\put(16201,-1936){\makebox(0,0)[lb]{\smash{{\SetFigFont{7}{8.4}{\familydefault}{\mddefault}{\updefault}$p$}}}}
\put(1501,-11311){\makebox(0,0)[lb]{\smash{{\SetFigFont{7}{8.4}{\familydefault}{\mddefault}{\updefault}$p$}}}}
\put(7651,-11311){\makebox(0,0)[lb]{\smash{{\SetFigFont{7}{8.4}{\familydefault}{\mddefault}{\updefault}$q$}}}}
\put(7651,-13036){\makebox(0,0)[lb]{\smash{{\SetFigFont{7}{8.4}{\familydefault}{\mddefault}{\updefault}$t'$}}}}
\put(7651,-15136){\makebox(0,0)[lb]{\smash{{\SetFigFont{7}{8.4}{\familydefault}{\mddefault}{\updefault}$t$}}}}
\put(4426,-16411){\makebox(0,0)[lb]{\smash{{\SetFigFont{7}{8.4}{\familydefault}{\mddefault}{\updefault}$r$}}}}
\put(1426,-16336){\makebox(0,0)[lb]{\smash{{\SetFigFont{7}{8.4}{\familydefault}{\mddefault}{\updefault}$u$}}}}
\put(2851,-12886){\makebox(0,0)[lb]{\smash{{\SetFigFont{7}{8.4}{\familydefault}{\mddefault}{\updefault}$R$}}}}
\put(17326,-11461){\makebox(0,0)[lb]{\smash{{\SetFigFont{7}{8.4}{\familydefault}{\mddefault}{\updefault}$q$}}}}
\put(17326,-15436){\makebox(0,0)[lb]{\smash{{\SetFigFont{7}{8.4}{\familydefault}{\mddefault}{\updefault}$t'$}}}}
\put(13276,-14761){\makebox(0,0)[lb]{\smash{{\SetFigFont{7}{8.4}{\familydefault}{\mddefault}{\updefault}$w$}}}}
\put(14251,-12661){\makebox(0,0)[lb]{\smash{{\SetFigFont{7}{8.4}{\familydefault}{\mddefault}{\updefault}$w'$}}}}
\put(1426,-20461){\makebox(0,0)[lb]{\smash{{\SetFigFont{7}{8.4}{\familydefault}{\mddefault}{\updefault}$p$}}}}
\put(1501,-25936){\makebox(0,0)[lb]{\smash{{\SetFigFont{7}{8.4}{\familydefault}{\mddefault}{\updefault}$u$}}}}
\put(4426,-25936){\makebox(0,0)[lb]{\smash{{\SetFigFont{7}{8.4}{\familydefault}{\mddefault}{\updefault}$r$}}}}
\put(7726,-24736){\makebox(0,0)[lb]{\smash{{\SetFigFont{7}{8.4}{\familydefault}{\mddefault}{\updefault}$t$}}}}
\put(3676,-25936){\makebox(0,0)[lb]{\smash{{\SetFigFont{7}{8.4}{\familydefault}{\mddefault}{\updefault}$r'$}}}}
\put(11026,-25861){\makebox(0,0)[lb]{\smash{{\SetFigFont{7}{8.4}{\familydefault}{\mddefault}{\updefault}$u$}}}}
\put(13801,-25936){\makebox(0,0)[lb]{\smash{{\SetFigFont{7}{8.4}{\familydefault}{\mddefault}{\updefault}$r$}}}}
\put(17326,-21136){\makebox(0,0)[lb]{\smash{{\SetFigFont{7}{8.4}{\familydefault}{\mddefault}{\updefault}$q$}}}}
\put(13051,-25936){\makebox(0,0)[lb]{\smash{{\SetFigFont{7}{8.4}{\familydefault}{\mddefault}{\updefault}$r'$}}}}
\put(10951,-24511){\makebox(0,0)[lb]{\smash{{\SetFigFont{7}{8.4}{\familydefault}{\mddefault}{\updefault}$f(u)$}}}}
\put(14251,-22186){\makebox(0,0)[lb]{\smash{{\SetFigFont{7}{8.4}{\familydefault}{\mddefault}{\updefault}$f(r')$}}}}
\put(16051,-21811){\makebox(0,0)[lb]{\smash{{\SetFigFont{7}{8.4}{\familydefault}{\mddefault}{\updefault}$f^2(r')$}}}}
\put(4726,-8161){\makebox(0,0)[lb]{\smash{{\SetFigFont{8}{9.6}{\familydefault}{\mddefault}{\updefault}{\sf{A}}}}}}
\put(4726,-17686){\makebox(0,0)[lb]{\smash{{\SetFigFont{8}{9.6}{\familydefault}{\mddefault}{\updefault}{\sf{C}}}}}}
\put(4726,-27361){\makebox(0,0)[lb]{\smash{{\SetFigFont{8}{9.6}{\familydefault}{\mddefault}{\updefault}{\sf{E}}}}}}
\put(14326,-27361){\makebox(0,0)[lb]{\smash{{\SetFigFont{8}{9.6}{\familydefault}{\mddefault}{\updefault}{\sf{F}}}}}}
\put(14326,-17761){\makebox(0,0)[lb]{\smash{{\SetFigFont{8}{9.6}{\familydefault}{\mddefault}{\updefault}{\sf{D}}}}}}
\put(14326,-8161){\makebox(0,0)[lb]{\smash{{\SetFigFont{8}{9.6}{\familydefault}{\mddefault}{\updefault}{\sf{B}}}}}}
\put(2776,-26536){\makebox(0,0)[lb]{\smash{{\SetFigFont{7}{8.4}{\familydefault}{\mddefault}{\updefault}$r''$}}}}
\put(7726,-20461){\makebox(0,0)[lb]{\smash{{\SetFigFont{7}{8.4}{\familydefault}{\mddefault}{\updefault}$q$}}}}
\put(10726,-22636){\makebox(0,0)[lb]{\smash{{\SetFigFont{7}{8.4}{\familydefault}{\mddefault}{\updefault}$f^2(u)$}}}}
\put(11026,-21061){\makebox(0,0)[lb]{\smash{{\SetFigFont{7}{8.4}{\familydefault}{\mddefault}{\updefault}$p$}}}}
\put(17326,-24736){\makebox(0,0)[lb]{\smash{{\SetFigFont{7}{8.4}{\familydefault}{\mddefault}{\updefault}$t$}}}}
\put(14401,-13711){\makebox(0,0)[lb]{\smash{{\SetFigFont{7}{8.4}{\familydefault}{\mddefault}{\updefault}$Q$}}}}
\end{picture}%

%% file: strip.pstex_t
\begin{picture}(0,0)%
\includegraphics{strip.pstex}%
\end{picture}%
\setlength{\unitlength}{1184sp}%
\begingroup\makeatletter\ifx\SetFigFont\undefined%
\gdef\SetFigFont#1#2#3#4#5{%
  \reset@font\fontsize{#1}{#2pt}%
  \fontfamily{#3}\fontseries{#4}\fontshape{#5}%
  \selectfont}%
\fi\endgroup%
\begin{picture}(16866,7965)(1168,-8293)
\put(6301,-3961){\makebox(0,0)[lb]{\smash{{\SetFigFont{7}{8.4}{\familydefault}{\mddefault}{\updefault}$w$}}}}
\put(1426,-6736){\makebox(0,0)[lb]{\smash{{\SetFigFont{7}{8.4}{\familydefault}{\mddefault}{\updefault}$u$}}}}
\put(1426,-1261){\makebox(0,0)[lb]{\smash{{\SetFigFont{7}{8.4}{\familydefault}{\mddefault}{\updefault}$p$}}}}
\put(2851,-6736){\makebox(0,0)[lb]{\smash{{\SetFigFont{7}{8.4}{\familydefault}{\mddefault}{\updefault}$r'$}}}}
\put(3976,-6736){\makebox(0,0)[lb]{\smash{{\SetFigFont{7}{8.4}{\familydefault}{\mddefault}{\updefault}$r$}}}}
\put(7651,-5236){\makebox(0,0)[lb]{\smash{{\SetFigFont{7}{8.4}{\familydefault}{\mddefault}{\updefault}$t$}}}}
\put(12301,-4036){\makebox(0,0)[lb]{\smash{{\SetFigFont{7}{8.4}{\familydefault}{\mddefault}{\updefault}$p$}}}}
\put(16201,-5461){\makebox(0,0)[lb]{\smash{{\SetFigFont{7}{8.4}{\familydefault}{\mddefault}{\updefault}$q$}}}}
\put(11776,-7036){\makebox(0,0)[lb]{\smash{{\SetFigFont{7}{8.4}{\familydefault}{\mddefault}{\updefault}$r$}}}}
\put(11026,-5911){\makebox(0,0)[lb]{\smash{{\SetFigFont{7}{8.4}{\familydefault}{\mddefault}{\updefault}$u$}}}}
\put(14476,-6361){\makebox(0,0)[lb]{\smash{{\SetFigFont{7}{8.4}{\familydefault}{\mddefault}{\updefault}$t$}}}}
\put(11326,-6511){\makebox(0,0)[lb]{\smash{{\SetFigFont{7}{8.4}{\familydefault}{\mddefault}{\updefault}$r'$}}}}
\put(4726,-8161){\makebox(0,0)[lb]{\smash{{\SetFigFont{8}{9.6}{\familydefault}{\mddefault}{\updefault}{\sf{A}}}}}}
\put(14326,-8161){\makebox(0,0)[lb]{\smash{{\SetFigFont{8}{9.6}{\familydefault}{\mddefault}{\updefault}{\sf{B}}}}}}
\put(3826,-4111){\makebox(0,0)[lb]{\smash{{\SetFigFont{7}{8.4}{\familydefault}{\mddefault}{\updefault}$w'$}}}}
\put(7726,-1261){\makebox(0,0)[lb]{\smash{{\SetFigFont{7}{8.4}{\familydefault}{\mddefault}{\updefault}$q$}}}}
\end{picture}%